\documentclass[11pt,reqno]{amsart}
\usepackage[catalan,spanish,english]{babel}
\usepackage[psamsfonts]{amssymb}
\usepackage{cmmib57}
\usepackage{exscale}
\usepackage{amsmath}
\usepackage{amscd}
\usepackage{latexsym}
\newtheorem{stat}{Statement}[section]
\newtheorem{thm}[stat]{Theorem}

\newtheorem{prop}[stat]{Proposition}
\newtheorem{lemma}[stat]{Lemma}

\numberwithin{equation}{section}

\newcommand {\IH}{\mbox{${\mathcal H}$}}

\newcommand{\beqn}{\begin{equation*}}
\newcommand{\eeqn}{\end{equation*}}
\newcommand{\beq}{\begin{equation}}
\newcommand{\eeq}{\end{equation}}

\newcommand{\half}{\frac{1}{2}}
\newcommand{\onehalf}{\frac{3}{2}}
\newcommand{\thalf}{\frac{t}{2}}

\newcommand{\ffi}{\varphi}
\newcommand{\dkm}{\Delta_k^m}
\allowdisplaybreaks
\def\1{{\rm 1 }\hskip -0.21truecm 1}
\begin{document}

\title[Approximation of rough paths of fractional Brownian motion]{ Approximation of rough paths\\
of fractional Brownian motion}

\author[A. Millet]
{Annie Millet $^{*,\dag}$}
\thanks{
* Supported by the grants BFM 2003-01345, HF 2003-006, Direcci\'on
General de Investigaci\'on, Ministerio de Educaci\'on y Ciencia, Spain.}
\address{Laboratoire de Probabilit\'es et Mod\`eles Al\'eatoires
(CNRS UMR 7599), Universit\'es Paris~6-Paris~7, Boite Courrier 188, 4 place Jussieu, 75252 Paris Cedex 05,
{\it  and } SAMOS-MATISSE,
Universit\'e Paris 1, 90 Rue de Tolbiac, 75634 Paris Cedex 13}
\email{amil@ccr.jussieu.fr}

\author[M. Sanz-Sol\'e]{Marta Sanz-Sol\'e $^{*}$}
\thanks{\dag 
Partially supported by the program SAB 2003-0082,
Direcci\'on General de
Universidades, Ministerio de Educaci\'on y Ciencia, Spain.}
\address
{Facultat de Matem\`atiques, Universitat de Barcelona, Gran Via 585, E-08007 Barcelona}
\email{marta.sanz@ub.edu}



\begin{abstract}
We consider a geometric rough path associated with a fractional Brownian motion with Hurst parameter
$H\in]\frac{1}{4}, \frac{1}{2}[$. We give an approximation result in a modulus type distance, up to the
second order, by means
of a sequence of rough paths lying above elements of the reproducing kernel Hilbert space.
\end{abstract}
\maketitle


\section{Introduction}
\label{S0}
Consider a $d$--dimensional fractional Brownian motion $W^H $ with Hurst parameter $H\in]\frac{1}{4},\half[\cup]\half,1[$
and integral representation
\beq
\label{a1}
W_t^H = \int_0^1 K^H(t,s)\, dB_s,
\eeq
where $K^H(t,s)=0$, if $s\ge t$ and for $0<s<t$,
\beq
\label{a2}
K^H(t,s) = c_H\, (t-s)^{H-\half} + s^{H-\half}F_1\left(\frac{t}{s}\right)
\eeq
with
\beq
\label{a3}
F_1(z) = c_H\left(\half-H\right)\int_0^{z-1} u^{H-\frac{3}{2}}\left(1-(u+1)^{H-\half}\right)\, du,
\eeq
for $z>1$ (see for instance \cite{amn}, equation (42)). In (\ref{a1}), $B$ denotes a standard $d$--dimensional
Brownian motion and in (\ref{a2}), (\ref{a3}), $c_H$ denotes a positive real constant depending on $H$.

Let $p\in ]1,4[$ be such that $pH>1$. In \cite{coutinquian}, it is proved that the sequence
of smooth rough paths based on linear interpolations of $W^H$ converges in the $p$--variation distance. The limit
defines a geometric rough path with roughness $p$ lying above $W^H$. We will call  this object
{\it the enhanced fractional Brownian motion}.

In the recent papers \cite{frizvictoir1}, \cite{coutinfrizvictoir}, the $p$--variation distance on rough paths
 is replaced by a strictly stronger,  {\it modulus type} distance defined as follows:
\beqn
\bar{d}_p(x,y)= \sup_{0\le s<t\le 1}
\left(\sum_{i=1}^{[p]} \frac{|x^{(i)}_{s,t}-y^{(i)}_{s,t}|}{(t-s)^{\frac{i}{p}}}
\right).
\eeqn
In \cite{coutinfrizvictoir}, it is proved that the enhanced fractional Brownian motion can actually be obtained by means
 of the $\bar{d}_p$ distance and
also that linear interpolations of $W^H$ define stochastic processes with values in $\mathcal{H}^H$,
the reproducing kernel Hilbert space
associated with $W^H$ (see Theorem 3.3 in \cite {du} for a description of this space).
Then, the authors state a characterization of the
topological support of the enhanced fractional Brownian motion among other results.

Our aim in this work is to give a new approximation of the enhanced fractional Brownian motion by means of a sequence
of geometric rough
paths which, unlike those based on linear interpolations, are not smooth, but also belong to $\mathcal{H}^H$.
For the sake of simplicity, we restrict
to $[p]=2$. We are pretty confident that our results extend to $[p]=3$; however, dealing with higher generality
would most likely produce a very technical paper. Our result, as is stated in Theorem \ref{t1.1},
provides in particular a new approximation of the L\'evy area of the fractional Brownian motion.

For any $m\in\mathbb{N}$, we consider the dyadic grid
 $\left(t_k^m=k2^{-m}, k=0,1,\dots,2^m\right)$ and set
 $\Delta_k^m=]t_{k-1}^m,t_k^m]$ and $\Delta_k^m B = B_{t_k^m}- B_{t_{k-1}^m}$.
Define $B(m)_0=0$ and for $t\in\Delta_k^m$,
$B(m)_t=B_{t_{k-1}^m} + 2^m(t-t_{k-1}^m) \Delta_k^m\,B$.
Our  approximation sequence is defined by
\beq
\label{1.1}
W(m)^H_t = \int_0^t K^H(t,s)\dot B(m)_s\,ds,
\eeq
$m\in\mathbb{N}$, where $\dot B(m)_s$ denotes the derivative with respect to $s$ of the path $s\mapsto B(m)_s$.
Notice that $W(m)^H\in \mathcal{H}^H$.

Let $K_m^H$ be the
orthogonal projection  of $K^H(t,\cdot)$ on the $\sigma$-field generated by
$(\Delta_k^m, k=1,\cdots,m)$. That is, for any $0<s<t\le 1$,
\beq
\label{km}
K_m^H(t,s) = \sum_{k=1}^{2^m}\,2^m\left(\int_{\Delta_k^m\cap]0,t]}K^H(t,u)\, du\right) \1_{\Delta_k^m}(s).
\eeq
We clearly have
\beq
\label{1.2}
W(m)^H_t = \int_0^1  K^H_m(t,s)\,dB_s.
\eeq
For $H\in ]\frac{1}{2},1[$,  we set ${\bf W}=({\bf W}_{s,t}= ( W_{s,t}^{(1)},
0\le s\le t\le 1)$,  ${\bf W(m)}= ({\bf W(m)}_{s,t}= ( W(m)_{s,t}^{(1)}, 0\le s\le t\le 1)$, while for
$H\in ]\frac{1}{4},\frac{1}{2}[$ we set ${\bf W}= ({\bf W}_{s,t}=
 ( W_{s,t}^{(1)}, W_{s,t}^{(2)}, \linebreak 0\le s\le t\le 1)$ and
  ${\bf W(m)}= ({\bf W(m)}_{s,t}= ( W(m)_{s,t}^{(1)}, 
   W(m)_{s,t}^{(2)}, 0\le s\le t\le 1)$,
 $m\ge 1$.
 \smallskip

 The main result of the paper states the convergence of ${\bf W(m)}$ to ${\bf W}$ in the $\bar{d}_p$-- metric for
 $p\in]1,3[$.
 For $p\in]1,2[$, the result is an almost trivial consequence of Lemma \ref{la2} which establishes  H\"older continuity
 in the $L^2[0,1]$ norm of the kernels $K^H$, $K_m^H$, respectively, and a control of the quadratic mean error
 in the approximation of $K^H$ by $K_m^H$. For $p\in[2,3[$, the approximation of the L\'evy area relies on representation
 formulas for the second order multiple integrals by means of the operator $K^*$ given in (\ref{2.1}) and introduced in
 \cite{amn} (see also \cite{coutinquian}). There are two fundamental ingredients. Firstly, Proposition \ref{p2}, giving
 the rate of convergence of the approximation
 at the second order level in the $L^q(\Omega)$--modulus norm; secondly, Lemma \ref{GRR}, an extension of the
 Garsia-Rademich-Rumsey Lemma for geometric rough paths of any roughness $p$. Other technical results used in the proofs,
  mainly on the kernels  $K^H$ and $K_m^H$,  are
 given in the Appendix.

For simplicity, in general we shall not write explicitly the dependence on $H$; thus $W$ stands for
$W^H$, $K(t,s)$ for $K^H(t,s)$, etc.
For any $q\in[1,\infty[$,
we denote by $\Vert\cdot\Vert_q$ the $L^q(\Omega)$--norm. We make
the convention $\sum_{k=a}^b x_k=0$ if $b<a$ and denote by $C$ positive constants
with possibly different values. For additional notions and notation on rough paths, we refer the reader to \cite{lyons}.

\section{Approximation result}
\label{S1}

For $p\in ]1,+\infty[$ we set $\tilde{d}_p=\bar{d}_{p\wedge 2}$, that is
\beqn
\tilde{d}_p(x,y)= \sup_{0\le s<t\le 1}
\left(\sum_{i=1}^{[p]\wedge 2} \frac{|x^{(i)}_{s,t}-y^{(i)}_{s,t}|}{(t-s)^{\frac{i}{p}}}
\right).
\eeqn
The purpose of this section is to prove the following approximation result.
\smallskip

\begin{thm}
\label{t1.1}
Let $H\in ]\frac{1}{4}, \frac{1}{2}[$,   $p\in]2,4[$ (resp. $H\in ]\frac{1}{2}, 1 [$,  $p\in ]1,2[$),
be such that $pH>1$  and $q\in [1,+\infty[$. The sequence
 $\big( \tilde{d}_p\left({\bf W(m)} , {\bf W}\right), m\ge 1\big)$, converges to 0 in
 $L^q(\Omega)$ and a.s. Thus for $H\in ]\frac{1}{2},1[$ and $p\in ]1,2[$,  if
${\mathcal G_p}$ denotes the set of
dyadic geometric rough paths endowed with the norm $\tilde{d}_p(0,.)$ and $P^H$ denotes
the law of the fractional Brownian motion $W^H$, then the triple $(X,{\mathcal H}^H,P^H)$
is an abstract Wiener space.
\end{thm}

The next Proposition provides the auxiliary result to state the approximation of the first
 component of the enhanced  fractional Brownian motion.

\begin{prop}
\label{p1} Let $0\leq s<t\leq 1$, $q\in[1,\infty[$.\\
\indent  (i)  For any $H\in]0, \half[$, $\lambda\in[0,H[$,
\beq
\label{1.3}
\left\Vert  W_{s,t}^{(1)}-W(m)^{(1)}_{s,t}\right\Vert_q \le C 2^{-m\lambda}
 \vert t-s\vert^{H-\lambda}.
 \eeq
\indent (ii) For any $H\in]\half,1[$,  $\varepsilon\in[0,H[$, $\mu\in ]0,\frac{\varepsilon}{H(2H+1)}[$,
\beq
\label{1.4}
\left\Vert  W_{s,t}^{(1)}-W(m)^{(1)}_{s,t}\right\Vert_q \le C  2^{-m \mu}
\vert t-s\vert^{H-\varepsilon}.
\eeq
\end{prop}
\begin{proof} By the hypercontractivity inequality, it suffices to prove the results for
 $q=2$. In this case, it is an easy consequence of the identity
\beqn
E\left(\left\vert W_{s,t}^{(1)}-W(m)^{(1)}_{s,t}\right\vert^2\right)=
\int_0^1 \left\vert\left(K(t,u)-K(s,u)\right)-\left(K_m(t,u)-K_m(s,u)\right)\right\vert^2 du
\eeqn
and of Lemma \ref{la2}. Indeed,  by (\ref{a16}), we have
\beqn
E\left(\left\vert W_{s,t}^{(1)}-W(m)^{(1)}_{s,t}\right\vert^2\right)\le C |t-s|^{2H}.
\eeqn
Hence, if  $t-s<2^{-m}$, we easily obtain (\ref{1.3}) and (\ref{1.4}).

Assume now  $H\in]0, \half[$ and $t-s\ge 2^{-m}$. By  (\ref{a17}), for $\epsilon \in [0,H]$,
 \beqn
E\left(\left\vert W_{s,t}^{(1)}-W(m)^{(1)}_{s,t}\right\vert^2\right)\le C 2^{-2mH}
\le C 2^{-2 m\epsilon} \vert t-s\vert^{2(H-\epsilon)}.
\eeqn
Hence, (\ref{1.3}) follows.

Let $H\in]\half,1[$ and  $t-s\ge 2^{-m}$.
Let $\alpha \in ]0,1[$~; then
(\ref{a16}) and (\ref{a18}) imply
\beqn
\left\Vert  W_{s,t}^{(1)}-W(m)^{(1)}_{s,t}\right\Vert_q \le C |t-s|^{H(1-\alpha)} 2^{-m\lambda\alpha},
\eeqn
with $\lambda\in]0,\frac{1}{2H+1}[$. By taking $\alpha=\frac{\varepsilon}{H}$, we obtain (\ref{1.4}) with
$\mu=\lambda \frac{\varepsilon}{H}$.
\end{proof}
\smallskip

Throughout  the rest of this section, $H\in]\frac{1}{4}, \half[$.
Following \cite{amn}, let ${\IH}_K$ denote the set of functions
$\ffi:[0,1]\to \mathbb R$ such that
\[ ||\ffi||_{K}^2=\int_0^1 \ffi(s)^2K(1,s)^2 ds +\int_0^1 ds\left(\int_s^1
|\ffi(t)-\ffi(s)|\,|K|(dt,s)\right)^2<+\infty .\]
For any $\ffi\in{\IH}_K$, $0<s<t$, set
\begin{align}
&K^{\ast}\left({\1}_{]s,t]}(\cdot)\left( \ffi_{\cdot}-\ffi_s\right)\right)(u)={\1}_{]0,s]}(u)\,
\int_s^t\left(\ffi_r-\ffi_s\right) K(dr,u)\nonumber\\
&\qquad +{\1}_{]s,t]}(u)\,\left(K(t,u)\left(\ffi_u-\ffi_s\right)+\int_u^t \left(\ffi_r-\ffi_u\right) K(dr,u)\right).
\label{2.1}
\end{align}
Following  \cite{coutinquian},
\beq
\label{2.2}
W_{s,t}^{(2)} = \int_0^1 K^{\ast}\left(\1_{]s,t]}(\cdot)\left(W_{\cdot}-W_s\right)\right)(u) dB_u +
\half \vert t-s\vert^{2H}.
\eeq
Moreover, by Theorem 9 in \cite{mss}, for $W(m)$ defined in (\ref{1.1}) we have
\beq
\label{2.3}
W(m)^{(2)}_{s,t} =  \int_0^1 K^{\ast}\left(\1_{]s,t]}(\cdot)\left(W(m)_{\cdot}-W(m)_s\right)
\right)(u) \dot{B}(m)_u du.
\eeq

\begin{prop}
\label{p2}
For each $m\in \mathbb{N}$,  $0<s<t\le 1$, $q\in[1,\infty[$,
\beq
\label{2.4}
\Vert W_{s,t}^{(2)}-W(m)^{(2)}_{s,t}\Vert_q \le C 2^{-m\mu} |t-s|^{2H-\varepsilon},
\eeq
for some positive constants $C$ and any $\varepsilon \in ]0, 2H - \frac{1}{2} [$ and
$\mu\in]0, \frac{\varepsilon}{2}[$.
\end{prop}

Before proving  this proposition, we give an equivalent expression for $W(m)^2_{s,t}$,
as follows.
The integration by parts formula of Malliavin calculus (see e.g. \cite{nualart},
 Equation (1.49)) and  (\ref{1.2}) yield
$W(m)^{(2)}_{s,t}=A_{s,t}^1(m)+A_{s,t}^2(m)$, with
\begin{align}
A_{s,t}^1(m)& = \sum_{k=1}^{2^m} \int_0^1  du\,\1_{\dkm}(u) 2^m K^{*}\left(\1_{]s,t]}(\cdot)
 \int_{\dkm} dB_r \left(W(m)_{\cdot}-W(m)_s\right)\right)(u),\label{2.5}\\
A_{s,t}^2(m)& = \sum_{k=1}^{2^m} \int_0^1  du\,\1_{\dkm}(u) 2^m K^{*}\left(\1_{]s,t]}(\cdot)
 \int_{\dkm} dr \left(K_m(\cdot,r)-K_m(s,r)\right)\right)(u).\label{2.6}
\end{align}
By 
 definition, for $r\in\dkm$, $K_m(t,r)=2^m\int_{\dkm\cap]0,t]}K(t,u) du =
 2^mK(\1_{\dkm})(t)$. Since $h:=K\left(\1_{\dkm}\right)\in \IH_K$, the duality relation
given in \cite{mss}, equation (58) and Lemma \ref{la6} yield
 \begin{align*}
A_{s,t}^2(m)&= \sum_{k=1}^{2^m} \int_0^1  dr \,\1_{\dkm}(r) 2^{2m} \int_0^1 \, du
\,\1_{\dkm}(u) K^{*}\left(\1_{]s,t]}(\cdot)K\left(\1_{\dkm}\right)_{s,\cdot}\right)(u)\\
& = \sum_{k=1}^{2^m} \int_0^1  dr\,\1_{\dkm}(r) 2^{2m} \int_s^t K\left(\1_{\dkm}\right)(du)
\left(K\left(\1_{\dkm}\right)(u)-K\left(\1_{\dkm}\right)(s)\right)\\
&= \sum_{k=1}^{2^m} \int_0^1  dr  \1_{\dkm}(r) 2^{2m}\frac{\left(K\left(\1_{\dkm}\right)(t)
-K\left(\1_{\dkm}\right)(s)\right)^2}{2}\\
& =\half\int_0^1 \,dr\,\vert K_m(t,r)-K_m(s,r)\vert^2 =
\frac{1}{2} \|W(m)^{(1)}_{s,t}\|_2^2.
\end{align*}
Thus, since $E |W_t-W_s|^2 = |t-s|^{2H}$, Schwarz's inequality, 
 (\ref{a16}), (\ref{a17}) imply
\beq
\label{ados}
\left\vert A_{s,t}^2(m)-\half \vert t-s\vert^{2H}\right\vert \le C 2^{-m\varepsilon}
\vert t-s\vert^{2H-\varepsilon},
\eeq
for some positive constant $C$ and $\varepsilon \in ]0,H[$.

 Hence, in order to establish (\ref{2.4}) it suffices to prove that for
any small parameter $\varepsilon \in ]0, 4H-1[$ and $\mu \in ]0,\varepsilon[$,
\beq
\label{2.7}
E\left(\left\vert\int_0^1 K^{\ast}\left(\1_{]s,t]}(\cdot)\left(W_{\cdot}-
W_s\right)\right)(u) dB_u - A_{s,t}^1(m)\right\vert^2\right) \le
 C 2^{- m{\mu}} \vert t-s\vert^{4H-\epsilon}.
\eeq
for all $m\geq 1$.
We devote the next lemmas to the proof of this convergence, 
using the  expression of   the operator
$K^{*}$ given in (\ref{2.1}).
\begin{lemma}
\label{l1.1}
For any $0\leq s<t \leq 1$, $m\ge 1$, we set
\begin{align*}
T_1(s,t)&= \int_0^s dB_u\, \left(\int_s^t (W_r-W_s)\, K(dr,u)\right),\\
T_1(s,t,m)&= \sum_{k=1}^{2^m}\int_{\dkm} dB_r \, 2^m \left(\int_{\dkm\cap]0,s]}  du\,
 \left(\int_s^t (W(m)_v-W(m)_s)\, K(dv,u)\right)\right).
\end{align*}
 Then for any $\epsilon\in]0,2H[$  and $\mu \in ]0, \epsilon [$,
 there exists   $C>0$ such that 
\beq
\label{2.8}
E\left(\left\vert T_1(s,t,m)-T_1(s,t)\right\vert^2\right)\le C 2^{-m \mu }
\vert t-s\vert^{4H-\epsilon}.
\eeq
\end{lemma}
\begin{proof}  Assume $s\in \Delta_I^m$, $I\ge 1$; we consider the  decomposition
$$E\left(\left\vert T_1(s,t,m)-T_1(s,t)\right\vert^2\right)\le C \sum_{j=1}^3
\tau_{1,j}(s,t,m),$$
with
\begin{align}
\tau_{1,1}(s,t,m)
&=\sum_{k\in\{1,I-1,I\}} E\Bigg(\Bigg\vert \int_{\dkm} dB_r \, 2^m\nonumber\\
&\quad \times \left(\int_{\dkm\cap]0,s]}  du\, \left(\int_s^t (W(m)_v-W(m)_s)\, K(dv,u)
\right)\right)\Bigg\vert^2\Bigg),\label{2.9}\\
\tau_{1,2}(s,t,m)&= \sum_{k\in\{1,I-1,I\}} E\left(\left\vert \int_{\dkm\cap]0,s]} dB_r
 \left(\int_s^t (W_v-W_s)\, K(dv,r)\right)\right\vert^2\right),\label{2.10}\\
\tau_{1,3}(s,t,m&)= E\Bigg(\Bigg\vert \sum_{k=2}^{I-2} \int_{\dkm} dB_r \, 2^m
 \int_{\dkm} du\,\Big(\int_s^t (W(m)_v-W(m)_s)\, K(dv,u)\nonumber\\
&-\int_s^t(W_v-W_s)K(dv,r)\Big)\Bigg\vert^2\Bigg).\label{2.11}
\end{align}

By  Lemma \ref{la7}, (\ref{a7}), Schwarz's inequality and (\ref{a16}),
any term in the right hand-side of (\ref{2.9}) is bounded as follows.
 Let $\varepsilon\in]0,2H[$, $\lambda \in ]\frac{1-(2H-\varepsilon)}{2}, \frac{1}{2}[$;
then  $2H-3+2\lambda <-1$, $1-2\lambda-(2H-\varepsilon)<0$ and
\begin{align*}
&E\left(\left\vert \int_{\dkm} dB_r \, 2^m \left(\int_{\dkm\cap]0,s]}  du\,
 \left(\int_s^t (W(m)_v-W(m)_s)\, K(dv,u)\right)\right)\right\vert^2\right)\\
\quad&\le C \int_{\dkm} dr \int_0^1 d\rho \left\vert 2^m\int_{\dkm\cap]0,s]}
 du \int_s^t (K_m(v,\rho)-K_m(s,\rho))\, K(dv,u)\right\vert^2\\
\quad&\le C \int_{\dkm} dr \int_0^1 d\rho\, 2^m \int_{\dkm\cap]0,s]} du\,
\left(\int_s^t dv\,|v-u|^{2H-3+2\lambda}  \right)\\
&\quad \times \left(\int_s^t dv\,
|K_m(v,\rho)-K_m(s,\rho)|^2   |v-u|^{-2\lambda} \right)\\
 \quad&\le C\int_{\dkm\cap]0,s]} du\, (s-u)^{2H-2 +2\lambda  } |t-s|^{2H}\,
\left( |t-u|^{1-2\lambda} - |s-u|^{1-2\lambda}\right)  \\
 \quad&\le C\int_{\dkm\cap]0,s]} du\, (s-u)^{2H-2 +2\lambda  } |t-s|^{2H}\,
|t-s|^{2H-\varepsilon}\, |s-u|^{1-2\lambda - (2H-\varepsilon)}  \\
 \quad&\le C |t-s|^{4H-\varepsilon} \int_{\dkm\cap]0,s]} du\, |s-u|^{\varepsilon -1}   \leq
 C 2^{-m\varepsilon } |t-s|^{4H-\varepsilon}.
\end{align*}

Each term of the right hand-side of (\ref{2.10}) can be studied using a similar strategy.
 Thus we obtain for $\varepsilon \in ]0,2H[$~:
\beq
\label{tau}
\tau_{1,1}(s,t,m)+\tau_{1,2}(s,t,m)\le C 2^{-m\varepsilon} |t-s|^{4H-\varepsilon}.
\eeq
Set for $s\geq 3\cdot  2^{-m}$, and hence $I\geq 4$,
\begin{align*}
X_r&=\sum_{k=2}^{I-2} \1_{\dkm}(r) 2^m \int_{\dkm} du\,
 \Big(\int_s^t (W(m)_v-W(m)_s)\, K(dv,u)\\
& - \int_s^t (W_v-W_s) K(dv,r)\Big) .
\end{align*}
Notice that $X_r=\int_0^1 g(r,\rho)\, dB_{\rho}$, with
\begin{align*}
g(r,\rho)&= \sum_{k=2}^{I-2} \1_{\dkm}(r) 2^m \int_{\dkm} du\,
 \Big(\int_s^t K(dv,u) \left(K_m(v,\rho)-K_m(s,\rho)\right)\\
& - \int_s^t K(dv,r)
\left(K(v,\rho)-K(s,\rho)\right) \Big).
\end{align*}
Hence, by Lemma \ref{la7} and Schwarz's inequality,
$ \tau_{1,3}(s,t,m)\le C( \tau_{1,3,1}(s,t,m) + \tau_{1,3,2}(s,t,m))$, with
\begin{align*}
\tau_{1,3,1}(s,t,m)&=\sum_{k=2}^{I-2}  2^m \int_{\dkm}dr\,\int_{\dkm}du\,\int_0^1 d\rho \\
&\quad\times\Big\vert\int_s^t\big(K_m(v,\rho)-K_m(s,\rho)\big)
\big(K(dv,u)-K(dv,r)\big)\Big\vert^2,
\end{align*}
\begin{align*}
\tau_{1,3,2}(s,t,m)&=\sum_{k=2}^{I-2}  \int_{\dkm} dr\,\int_0^1 d\rho \Big\vert\int_s^t {K(dv,r)}\\
 &\quad\times \big(K_m(v,\rho)-K_m(s,\rho)-K(v,\rho)+K(s,\rho)\big) \Big\vert^2.
\end{align*}
Owing to (\ref{a7}), (\ref{a10}), we have for $\lambda\in]0,1[$, $u,r\in\dkm$,
\begin{align}
&\left\vert \frac{\partial K}{\partial v}(v,u)-\frac{\partial K}{\partial v}(v,r)
\right\vert\nonumber\\
&\quad\le C\left\vert \frac{\partial K}{\partial v}(v,u)-\frac{\partial K}{\partial v}(v,r)
\right\vert^\lambda
\left( \left\vert \frac{\partial K}{\partial v}(v,u)\right\vert^{1-\lambda}
+\left\vert \frac{\partial K}{\partial v}(v,r)
\right\vert^{1-\lambda}\right)\nonumber\\
&\quad\le C  2^{-m\lambda}|v-(u\vee r)|^{H-\onehalf} \left[ (u\wedge r)^{-1}
+|v-(u\vee r)|^{-1}\right]^{\lambda} .
\label{bound}
\end{align}

\noindent Thus, taking $\lambda:=H$ yields $\tau_{1,3,1}(s,t,m)\le C
 2^{-2mH} \sum_{j=1}^2 \tau_{1,3,1,j}(s,t,m)$, with
\begin{align*}
\tau_{1,3,1,1}(s,t,m)&=\sum_{k=2}^{I-2}2^m \int_{\dkm}dr \int_{\dkm}du \int_0^1 d\rho
\Big(\int_s^t dv |K_m(v,\rho)-K_m(s, \rho )| \\
&\quad \times |v-(u\vee r)|^{H-\onehalf}(u\wedge r)^{-H}\Big)^2, \\
\tau_{1,3,1,2}(s,t,m)&=\sum_{k=2}^{I-2}2^m \int_{\dkm}dr \int_{\dkm}du \int_0^1 d\rho
\Big(\int_s^t dv |K_m(v,\rho)-K_m(s,\rho)|\\
&\quad \times |v-(u\vee r)|^{-\onehalf}\Big)^2.
\end{align*}
Let $a=2-\epsilon$, with $\epsilon\in]0,2H[$. Schwarz's inequality
along with  (\ref{a16}) yield
\begin{align}
\tau_{1,3,1,1}(s,t,m)&\le C\sum_{k=2}^{I-2}2^m \int_{\dkm} dr \int_{\dkm} du
\left(\int_s^t dv |v-(u\vee r)|^{-a} dv \right)\nonumber\\
&\quad\times \left(\int_s^t dv |v-s|^{4H-3+a}|u\wedge r|^{-2H}\right)\nonumber\\
&\le C |t-s|^{4H-\epsilon}\int_{t_1^m}^{t_{I-2}^m} du (s-\overline{u}_m)^{\epsilon-1}
 (\underline{u}_m)^{-2H} \nonumber \\
&\le C |t-s|^{4H-\epsilon} \,  s^{\epsilon -2H}
\le |t-s|^{4H-\epsilon}\, 2^{-m(\epsilon -2H)}.\label{2.12}
\end{align}
Indeed, $\int_s^t dv |v-(u\vee r)|^{-2+\epsilon} \le C(s-\overline{u}_m)^{\epsilon-1}$
for $\overline{u}_m$ defined by (\ref{numero}).
Let  $\epsilon \in ]0, 2H[$
 using Schwarz's  inequality and (\ref{a16}), we obtain
\begin{align}
\tau_{1,3,1,2}(s,t,m)&\le C \sum_{k=2}^{I-2}2^m \int_{\dkm}dr \int_{\dkm}du
 \Big(\int_s^t dv |v-(u\vee r)|^{-2-2H+\varepsilon}\Big)\nonumber\\
&\quad\times  \Big(\int_s^t dv |v-(u\vee r)|^{2H-\varepsilon-1}|v-s|^{2H }\Big)\nonumber\\
&\le  C |t-s|^{4H-\varepsilon}\int_{t_1^m}^{t_{I-2}^m} du \int_s^t  dv
 (v-\overline{u}_m)^{-2-2H+\varepsilon}\nonumber\\
&\le C  |t-s|^{4H-\varepsilon}\int_{t_1^m}^{t_{I-2}^m} du (s-\overline{u}_m)^{-1-2H+\varepsilon}\nonumber\\
&\le C  |t-s|^{4H-\varepsilon} 2^{-m(\varepsilon-2H)}.\label{2.13}
\end{align}

From (\ref{2.12}), (\ref{2.13}) we deduce that for $\epsilon \in ]0, 2H[$,
\beq
\label{2.14}
\tau_{1,3,1}(s,t,m)\le C |t-s|^{4H-\epsilon}\, 2^{-m\epsilon}.
\eeq
Let $\delta \in ]0,2H[$, $\alpha\in ]0,2H[$, $\lambda \in ]0,1[$ and $\mu\in ]\frac{1}{2},1-H[$.
 Notice that for these choices, $-2\mu+1-2H+\delta<0$.
H\"older's inequality together with (\ref{a16}) and (\ref{a17}) yield for any $\lambda\in]0,1[$,
\beqn
\tau_{1,3,2}(s,t,m)\le C \tau_{1,3,2,1}(s,t,m)^\lambda \, \tau_{1,3,2,2}(s,t,m)^{1-\lambda},
\eeqn
where
\begin{align*}
 \tau_{1,3,2,1}(s,t,m) &= \int_{t_1^m}^{t_{I-2}^m} dr \left(\int_s^t dv (v-r)^{2H-3+2\mu}\right)
 \left(\int_s^t dv (v-r)^{-2\mu} (v-s)^{2H}\right),\\
 \tau_{1,3,2,2}(s,t,m) &= \int_{t_1^m}^{t_{I-2}^m} dr \left(\int_s^t dv (v-r)^{2H-3+2\mu}\right)
 \left(\int_s^t dv (v-r)^{-2\mu} 2^{-2mH}\right).
 \end{align*}
For the first term we have
\begin{align*}
 \tau_{1,3,2,1}(s,t,m)&\le C\vert t-s\vert^{4H-\delta}\int_{t_1^m}^{t_{I-2}^m} dr
(s-r)^{2H-2+2\mu} (s-r)^{-2\mu+1-2H+\delta}\\
 &\le C \vert t-s\vert^{4H-\delta},
 \end{align*}
while for the second one, we obtain
\begin{equation*}
\tau_{1,3,2,2}(s,t,m) \le C 2^{- 2mH}\vert t-s\vert^{2H-\alpha}
 \int_{t_1^m}^{t_{I-2}^m} dr (s-r)^{2H-2+2\mu}(s-r)^{-2\mu+1-2H+\alpha}.
 \end{equation*}
Consequently,
\beqn
\tau_{1,3,2}(s,t,m)\le C\,  |t-s|^{(4H-\delta)\lambda + (2H-\alpha)(1-\lambda)} \,  2^{-2m H (1-\lambda)}\, .
\eeqn
Take $\alpha, \delta$ arbitrarily small and
 $1-\lambda = \frac{\epsilon-H\delta}{2H+\alpha}$. Then
for $\beta < \epsilon < 2H$,  we have proved that
$$\tau_{1,3,2}(s,t,m) \leq C |t-s|^{4H-\epsilon} 2^{-m\beta}.$$
 This inequality, together with (\ref{tau}) and (\ref{2.14}) yields
 (\ref{2.8}).
\end{proof} 

\begin{lemma}
\label{l1.2}
For any $0\leq s<t\leq 1$, set
\begin{align*}
T_2(s,t)&=\int_s^t dB_u\,K(t,u) (W_u-W_s)\, , \\
T_2(s,t,m)&= \sum_{k=1}^{2^m}\int_{\dkm} dB_r \, 2^m \left(\int_{\dkm\cap]s,t]}
 du K(t,u)\left( W(m)_u-W(m)_s\right)\right).
\end{align*}
Then,  for $b\in ]0,2H[$, there exists a constant $C>0$ such that for each $m\ge 1$
\beq
\label{221}
E\left(\left\vert T_2(s,t,m)-T_2(s,t)\right\vert^2\right)\le C 2^{-mb}\vert t-s\vert^{4H-b}.
\eeq
\end{lemma}
\begin{proof}  Let $s\in \Delta_I^m$, $t\in\Delta_J^m$.
We have
\beqn
E\left(\left\vert T_2(s,t,m)-T_2(s,t)\right\vert^2\right)\le C\sum_{j=1}^3 T_{2,j}(s,t,m),
\eeqn
with for ${\mathcal I}= \{I, I+1, J-2, J-1J\}$
\begin{align*}
T_{2,1}(s,t,m)&=\sum_{k\in {\mathcal I}} E\Big(\Big\vert\int_{\dkm}\! dB_r\, 2^m
\int_{\dkm\cap]s,t]} \! du K(t,u) 
\left( W(m)_u-W(m)_s\right)\Big\vert^2\Big),\\
T_{2,2}(s,t,m)&=\sum_{k\in {\mathcal I}}
 E\Big(\Big\vert\int_{\dkm\cap]s,t]} dB_r\, K(t,r)(W_r-W_s)\Big\vert^2\Big),\\
T_{2,3}(s,t,m)&=E\Big(\Big\vert \sum_{k=I+2}^{J-3} \int_{\dkm} dB_r
\Big[  \,2^m \int_{\dkm\cap]s,t]} du K(t,u)\left( W(m)_u-W(m)_s\right)\\
&\qquad -K(t,r)(W_r-W_s)\Big] \Big\vert^2\Big).
\end{align*}
Owing to Lemma \ref{la7} applied to the Gaussian process
$$X_r:= \1_{\dkm}(r) \int_0^1 dB_{\rho} \left(2^m\int_{\dkm\cap]s,t]} du
 K(t,u)\big(K_m(u,\rho)-K_m( s,\rho)\big)\right)$$
and Schwarz's  inequality,
we have for any $k=1,\cdots, 2^m$,
\begin{align*}
T(s,t,m,k)&:=E\left(\left\vert\int_{\dkm} dB_r\, 2^m \int_{\dkm\cap]s,t]} du K(t,u)
\left( W(m)_u-W(m)_s\right)\right\vert^2\right)\\
&\quad \le C 2^{2m} \int_{\dkm} dr \int_0^1 d\rho \left(\int_{\dkm\cap]s,t]} du K^2(t,u)\right)\\
&\qquad \times \left( \int_{\dkm\cap]s,t]} du\left\vert K_m(u,\rho)-K_m(s,\rho)\right\vert^2\right).
\end{align*}
Let $k=I,I+1$;
since $\int_{\dkm\cap]s,t]} du K^2(t,u)\le \int_{]s,t]} du K^2(t,u)\le C|t-s|^{2H}$,
we have for any $b\in]0,2H[$,
\begin{align*}
T(s,t,m,k)& \le  C 2^{m}  |t-s|^{2H} \left(\int_{\dkm\cap]s,t]}du |u-s|^{2H}\right)\\
&\le C 2^m |t-s|^{4H-b}\int_{\dkm\cap]s,t]}du |u-s|^b\le C 2^{-mb}|t-s|^{4H-b}.
\end{align*}
Let $k = J-2,J-1,J$  with $J-2>I+1$ then for $u\in\dkm$,
 (\ref{a8}) implies $|K(t,u)|^2\leq C |t-u|^{2H-1}$ and $|t-u|\le C 2^{-m}$;
 we obtain for $b\in]0,2H[$,
\begin{align*}
T(s,t,m,k)& \le  C 2^{m} \left(\int_{\dkm\cap]s,t]} du |t-u|^{2H-1-b} 2^{-mb} du \right)
 \left(\int_{\dkm\cap]s,t]} du |u-s|^{2H}\right)\\
&\le C |t-s|^{4H-b} 2^{-mb}.
\end{align*}
We therefore have proved that for $b\in ]0, 2H[$,
\beq
\label{222}
T_{2,1}(s,t,m)\le C 2^{-bm}|t-s|^{4H-b}.
\eeq

The analysis of the term $T_{2,2}(s,t,m)$ is easier. Indeed, the isometry property
of the stochastic integral yields
for any $k=1,\cdots, 2^m$,
\beq
\label{223}
E\left(\left\vert\int_{\dkm\cap]s,t]} dB_r\, K(t,r)(W_r-W_s)\right\vert^2\right)
 = C \int_{\dkm\cap]s,t]} dr\, K^2(t,r) |r-s|^{2H}.
\eeq
For the particular values of $k\in\mathcal{I}$,
 the right hand-side of (\ref{223})
can be analyzed  following similar ideas as for
$T_{2,1}(s,t,m)$, which yields for $b\in ]0, 2H[$
\beq
\label{224}
T_{2,2}(s,t,m) \le C  2^{-mb} |t-s|^{4H-b}.
\eeq
We now study $T_{2,3}(s,t,m)$ and note that $T_{2,3}(s,t,m)=0$ if $|t-s|\leq 2^{-m}$.
Thus, we may assume that
$t-s\geq 2^{-m}$.  First, we apply Lemma \ref{la7} and  obtain
\beqn
T_{2,3}(s,t,m)\le C(T_{2,3,1}(s,t,m)+T_{2,3,2}(s,t,m)),
\eeqn
where
\begin{align*}
T_{2,3,1}(s,t,m)&= \int_{\overline{s}_m}^{\underline{t}_m-2^{1-m}} dr \int_0^1 d\rho
\Big\vert  2^m\int_{\underline{r}_m}^{\overline{r}_m} du \Big(K(t,u)-K(t,r)\Big)\\
 &\qquad \times \Big(K_m(u,\rho)-K_m(s,\rho)\Big)\Big\vert^2,\\
T_{2,3,2}(s,t,m)&= \int_{\overline{s}_m}^{\underline{t}_m-2^{1-m}} dr \int_0^1 d\rho
 \Big\vert 2^m \int_{\underline{r}_m}^{\overline{r}_m} du K(t,r)\\
&\quad \times  \Big( \big[K_m(u,\rho)  -K_m(s,\rho)\big]-\big[K(r,\rho)-K(s,\rho)\big]\Big)
\Big\vert^2.
\end{align*}
By Schwarz's inequality and (\ref{a16}), for $b\in ]0,2H[$,
\begin{align*}
T_{2,3,1}(s,t,m)&\le \int_{\overline{s}_m}^{\underline{t}_m-2^{1-m }} dr
2^m\int_{\underline{r}_m}^{\overline{r}_m} du
\vert K(t,u)-K(t,r)\vert^2 \vert u-s\vert^{2H}\\
&\le C |t-s|^{2H} \int_{\overline{s}_m}^{\underline{t}_m-2^{1-m}} dr
2^m\int_{\underline{r}_m}^{\overline{r}_m} du
\vert K(t,u)-K(t,r)\vert^2\\
&\le C 2^{-2mH}  |t-s|^{2H} \leq C 2^{-mb} |t-s|^{4H-b}
\end{align*}
where the last inequalities  follow from (\ref{forlemma}) and $|t-s|\geq 2^{-m}$.

Owing to (\ref{a17}), we have for $u\in [\underline{r}_m,\overline{r}_m]$
\begin{align*}
 \int_0^1 d\rho \vert
K_m(s,\rho)-K(s,\rho)\vert^2 \le & C 2^{-2mH},\\
 \int_0^1 d\rho \vert
 K_m(u,\rho)-K(r,\rho)\vert^2 \le& C \int_0^1 d\rho
\Big( \vert  K_m(u,\rho)-K(u,\rho)\vert^2\\
&\quad +   \vert  K(u,\rho)-K(r,\rho)\vert^2\Big) \leq C 2^{-2mH}.
\end{align*}
Schwarz's inequality, along with (\ref{a8}) and the above  estimates yield
\begin{align*}
T_{2,3,2}(s,t,m)&\le C \int_{\overline{s}_m}^{\underline{t}_m-2^{1-m}} dr
 2^{-2mH} \left(|r|^{2H-1}+|t-r|^{2H-1}\right)\\
&\le C 2^{-2mH}\left(t^{2H}-s^{2H}+|t-s|^{2H}+2^{-2mH}\right)\\
&\le C 2^{-2mH}|t-s|^{2H} \leq C 2^{-mb} |t-s|^{4H-b}
\end{align*}
for $b\in ]0,2H[$. Indeed, for each $H\in]0,\half[$, and $s<t$, $t^{2H}-s^{2H}
\le (t-s)^{2H}$  and we are assuming that $2^{-m}<|t-s|$.
Thus, (\ref{221}) is proved.
\hfill\qed

\begin{lemma}
\label{l1.3}
For any $0\leq s<t\leq 1$, set
\begin{align*}
T_3(s,t)&=\int_s^t dB_u\,\int_u^t  K(dr,u)(W_r-W_u)\\
T_3(s,t,m)&= \sum_{k=1}^{2^m} 2^m  \int_{\dkm} dB_r \, \int_{\dkm\cap]s,t]}  du \int_u^t K(dv,u)
\left( W(m)_v-W(m)_u\right).
\end{align*}
There exists a positive constant $C$ such that, for any $\epsilon\in]0,4H-1[$
\beq
\label{2.15}
E\left(\left\vert T_3(s,t,m)-T_3(s,t)\right\vert^2\right)\le C 2^{-m\epsilon}\vert t-s\vert^{4H-\epsilon},
\eeq
for each $m\ge 1$.
\end{lemma}
\begin{proof} Assume $s\in\Delta_I^m$, $t\in \Delta_J^m$; we consider the upper bound
\beqn
E\left(\vert T_3(s,t,m)-T_3(s,t)|^2\right) \le C \sum_{j=1}^3 T_{3,j}(s,t,m),
\eeqn
where  for ${\mathcal J}= \{ I, I+1,J-1,J\}$
\begin{align}
T_{3,1}(s,t,m)&=\sum_{k\in {\mathcal J}}  E\Big(\big\vert 2^m \int_{\dkm} dB_r
 \int_{\dkm\cap]s,t]} du \int_u^t K(dv,u)\nonumber\\
&\quad \times (W(m)_v-W(m)_u)\big\vert^2 \Big),\label{2.16}\\
T_{3,2}(s,t,m)&=\sum_{k\in {\mathcal J}} E\Big(\Big\vert \int_{\dkm\cap]s,t]} dB_r
\int_r^t  K(dv,r)(W_v-W_r)\Big\vert^2\Big),\label{2.17}\\
T_{3,3}(s,t,m)&=E\Big(\Big\vert\sum_{k=I+2}^{J-2} 2^m \int_{\dkm} dB_r\int_{\dkm}
du\nonumber\\
&\quad\times \Big(\int_u^t K(dv,u)(W(m)_v-W(m)_u)-\int_r^t K(dv,r)(W_v-W_r)\Big)
\Big\vert^2\Big).\nonumber
\end{align}
Lemma \ref{la7} along with Schwarz's inequality yield for each term of the sum in the right
hand side of (\ref{2.16}) the upper bound
\beqn
C\int_{\dkm} dr \int_0^1 d\rho\, 2^m\int_{\dkm\cap]s,t]} du \left(\int_u^t K(dv,u)
\big(K_m(v,\rho)-K_m(u,\rho)\big)\right)^2.
\eeqn
Fix $a\in]2-4H,1]$. From Schwarz's inequality, (\ref{a7}) and (\ref{a16})
 we deduce the following estimates for this integral:
\begin{align*}
&C \int_{\dkm} dr 
2^m\int_{\dkm\cap]s,t]} du \left(\int_u^t dv  |v-u|^{-a}\right)
\left(\int_u^t |v-u|^{4H-3+a}\right)\\
&\quad \le C \left(2^{-m}\wedge |t-s|\right) |t-s|^{4H-1}.
\end{align*}
A similar analysis yields  the same result for each term in the right hand-side of
(\ref{2.17}). Consequently,
\beq
\label{2.18}
T_{3,1}(s,t,m)+T_{3,2}(s,t,m) \le C \left( 2^{-m}\wedge |t-s|\right) |t-s|^{4H-1}.
\eeq
If $|t-s|\leq 2^{-m}$ then $T_{3,3}(s,t,m)=0$. Hence, let us assume that $t-s\geq 2^{-m}$;
in this case $T_{3,3}(s,t,m)$ is equal to
 $E\left(\int_0^1 dB_r X_r\right)^2$, with
$X_r= \int_0^1 dB_\rho g(r,\rho)$, and
\begin{align*}
g(r,\rho)&=  \sum_{k=I+2}^{J-2} \1_{\dkm} (r)2^m \int_{\dkm} du \Big[\int_u^t
K(dv,u)(K_m(v,\rho)-K_m(u,\rho))\\
&\quad -\int_r^t K(dv,r)(K(v,\rho)-K(r,\rho))\Big].
\end{align*}
We at first study the contribution to  $T_{3,3}(s,t,m)$ of the integrands
\begin{align*}
g_1(r,\rho)&= \sum_{k=I+2}^{J-2} \1_{\dkm} (r) 2^m \int_{\dkm} du \int_u^{u\vee r}
K(dv,u)(K_m(v,\rho)-K_m(u,\rho)),\\
g_2(r,\rho)&=  \sum_{k=I+2}^{J-2} \1_{\dkm} (r) 2^m \int_{\dkm} du \int_r^{u\vee r}
K(dv,r)(K(v,\rho)-K(r,\rho)),
\end{align*}
which we denote by $T_{3,3,j}(s,t,m)$, $j=1,2$. Actually, both are similar and therefore
 we only study the first one.
Lemma \ref{la7}, (\ref{a7}), (\ref{a16}) and Schwarz's  inequality imply,
 for each $a\in]2-4H,1]$,
\begin{align}
T_{3,3,1}(s,t,m)&\le C \!  \sum_{k=I+2}^{J-2}  \!  2^m\int_{\dkm}\!  dr  \int_{\dkm}\! du
\int_u^{u\vee r}\! dv |v-u|^{-a} \int_u^{u\vee r }\! dv |v-u|^{4H-3+a}\nonumber\\
&\le C 2^{-m(4H-1)} |t-s|. \label{2.19}
\end{align}
We end the analysis of the term $T_{3,3}(s,t,m)$ by studying the contribution of
  $T_{3,3,3}((s,t,m)$ defined in terms of the integrand
\begin{align*}
g_3(r,\rho)&=  \sum_{k=I+2}^{J-2} \int_{\dkm} dr  2^m \int_{\dkm} du\int_{u\vee r}^t
\Big[  K(dv,u)(K_m(v,\rho)-K_m(u,\rho))\\
&\quad -K(dv,r)(K(v,\rho)-K(r,\rho))\Big].
\end{align*}
Notice that $g_3(r,\rho)$ is the sum of two analogous terms where the set $\dkm$
of the integral with respect to the variable $u$
is replaced by $
[\underline r_m,r[$, $
[r,\overline r_m[$, respectively.
Again, the contribution of both terms is similar,
so that we concentrate on the first one. That is, we consider
\begin{align*}
T_{3,3,3}^{+}(s,t,m):&=E\Big(\Big\vert  \sum_{k=I+2}^{J-2} 2^m \int_{\dkm} dB_r
 \int_{
[\underline{r}_m,r[} du  \int_r^t \Big[ K(dv,u)\\
&\quad \times (W(m)_v-W(m)_u)-K(dv,r)(W_v-W_r) \Big] \Big\vert^2\Big).
\end{align*}
As before, all the arguments rely on Lemma \ref{la7}, (\ref{a7}), (\ref{a16}),
a suitable factorization of the integrands along with Schwarz's  inequality.
In order to deal with the singularity at $v=r$, we first replace the integral with
respect to the variable $v$ by $\int_r^{\overline r_m+2^{-m}}$. Given $a\in ]2-4H,1[$,
 the corresponding  contribution to  $T_{3,3,3}^+(s,t,m)$ is bounded by
\begin{align}
&C \sum_{k=I+2}^{J-2}  2^m \int_{\dkm} dr \int_{
[\underline r_m,r[} du
\int_0^1 d\rho\Big( \Big| \int_r^{\overline r_m+2^{-m}} K(dv,u)\nonumber\\
&\quad\times (K_m(v,\rho)-K_m(u,\rho))\Big|^2 +\Big|\int_r^{\overline{r}_m+2^{-m}}
 K(dv,r)(K(v,\rho)-K(r,\rho))\Big|^2\Big)\nonumber\\
&\le C \sum_{k=I+2}^{J-2} 2^m \int_{\dkm} \!\! dr \int_{
[\underline r_m,r[} \!\!  du
\int_r^{\overline{r}_m+2^{-m}}\!\! dv |v-r|^{-a}
\int_r^{\overline{r}_m+2^{-m}}\!\!  dv |v-r|^{4H-3+a}\nonumber\\
&\le C 2^{-m(4H-1)}|t-s|. \label{2.20}
\end{align}

Let us finally consider the range $]r_m+2^{-m},t[$ for the variable $v$.
We have to study two terms:
\begin{align*}
M_1(s,t,m)&= \sum_{k=I+2}^{J-2} 2^m \int_{\dkm} dr \int_{
[\underline r_m,r[} du \int_0^1 d\rho \Big(\int_{\overline{r}_m+2^{-m}}^t dv\\
&\quad\times \vert K_m(v,\rho)-K_m(u,\rho)\vert
\left\vert \frac{\partial K}{\partial v}(v,u)-\frac{\partial K}{\partial v}(v,r)
\right\vert\Big)^2,\\
M_2(s,t,m)&=  \sum_{k=I+2}^{J-2} 2^m \int_{\dkm} dr \int_{
[\underline r_m,r[} du \int_0^1 d\rho \Big(\int_{\overline{r}_m+2^{-m}}^t dv
 \Big| \frac{\partial K}{\partial v}(v,r)\Big|\\
&\quad \times \big[(K_m(v,\rho)-K_m(u,\rho))-(K(v,\rho)-K(r,\rho))\big]\Big)^2.
\end{align*}

For $M_1(s,t,m)$, we proceed in a similar way as for the term $\tau_{1,3,1}(s,t,m)$
in Lemma \ref{l1.1},
as follows. By means of (\ref{bound}) we obtain for $\lambda\in]0,1[$
 $M_1(s,t,m)\le C 2^{-2m\lambda}\left(M_{1,1}(s,t,m)+M_{1.2}(s,t,m)\right)$, with
\begin{align*}
M_{1,1}(s,t,m)&= \sum_{k=I+2}^{J-2}  2^m \int_{\dkm} dr \int_{
[\underline{r}_m,r[} du
\, u^{-2\lambda}\int_0^1 d\rho \Big(\int_{\overline{r}_m+2^{-m}}^t dv\\
&\qquad \vert K_m(v,\rho)-K_m(u,\rho)\vert |v-r|^{H-\onehalf}\Big)^2,\\
M_{1,2}(s,t,m)&= \sum_{k=I+2}^{J-2} 2^m \int_{\dkm} dr \int_{
[\underline{r}_m,r[} du \int_0^1 d\rho \Big(\int_{\overline{r}_m+2^{-m}}^t dv\\
&\qquad \vert K_m(v,\rho)-K_m(u,\rho)|v-r|^{H-\onehalf-\lambda}\Big)^2.
\end{align*}
Let $a\in]2-4H,1[$, $\lambda\in]0,\half[$. Since $t-s\geq 2^{-m}$,  for
$u\in [\underline{r}_m,r[$, we have
\beqn
\int_{\overline r_m+2^{-m}}^t dv |v-r|^{2H-3+a}|v-u|^{2H} 
\leq C |t-r|^{4H+a-2}.
\eeqn
Consequently,  since $r\geq u\geq \underline{r}_m \geq t_{I+1}$ implies $u\geq \frac{r}{2}$
\begin{align}
M_{1,1}&(s,t,m)\le C \sum_{k=I+2}^{J-2} 2^m \int_{\dkm} dr
\int_{
[\underline{r}_m,r[} du\, u^{-2\lambda}\Big(\int_{\overline{r}_m+2^{-m}}^t dv
 |v-r|^{-a}\Big)\nonumber\\
&\qquad\times \Big( \int_{\overline{r}_m+2^{-m}}^t dv |v-r|^{2H-3+a}|v-u|^{2H}\Big)\nonumber\\
&  \leq C \int_s^{t} \!\! r^{-2\lambda} |t-s|^{4H-1} dr
 \le C |t-s|^{4H-2\lambda}.\label{2.21}
\end{align}
Analogously, for $b\in]2+2\lambda-4H,1[$, $\lambda\in]0,2H-\half[$
and $|t-s|\geq 2^{-m}$
\begin{align}
M_{1,2}(s,t,m)&\le C  \sum_{k=I+2}^{J-2} 2^m \int_{\dkm} dr
\int_{ [\underline{r}_m,r[} du \Big(\int_{\overline r_m+2^{-m}}^t dv
 |v-r|^{-b}\Big)\nonumber\\
&\quad\times \Big( \int_{\overline{r}_m+2^{-m}}^t dv |v-r|^{2H-3-2\lambda+b}
|v-u|^{2H}\Big)\nonumber\\
& \le C \int_s^t |t-r|^{4H-1-2\lambda} dr = C |t-s|^{4H-2\lambda}.
\label{2.22}
\end{align}
Finally, if we additionally use (\ref{a17}), we obtain for   $a\in ]2-4H,1[$
\begin{align}
M_2(s,t,m)&\le C  \sum_{k=I+2}^{J-2} 2^m \int_{\dkm} dr \int_{[\underline r_m,r[}
 du\Big(\int_r^t dv |v-r|^{-a}\Big)\nonumber\\
&\quad \times \Big(\int_{ \underline{r}_m +2^{-m}}^t dv |v-r|^{2H-3+a} 2^{-2mH}\Big)\nonumber\\
&\le C \int_s^t |t-r|^{1-a}\, 2^{-m(4H-2+a)} dr  \leq C 2^{-mb} |t-s|^{4H-b}  \label{2.23}
\end{align}
 for$b\in ]0,4H-1[$.  We easily check that (\ref{2.15}) follows from
 (\ref{2.18})--(\ref{2.23}).
\end{proof} 
\smallskip

\noindent{\it Proof of Proposition \ref{p2}}: We remark that Lemmas \ref{l1.1} to
 \ref{l1.3} yield the upper bound (\ref{2.7}).
Therefore, for $q=2$, (\ref{2.4}) follows from (\ref{ados}) and (\ref{2.7}).
 The hypercontractivity
inequality yields the validity of the same inequality for any
$q\in]2,\infty[$. \hfill $\square$
\smallskip

\noindent{\it Proof of Theorem \ref{t1.1}}:

Let $H\in ]\frac{1}{2},1[$ and $p\in ]\frac{1}{H},2[$.
The convergence of $\tilde{d}_p(\bf {W(m)},\bf W)$ to zero in $L^q(\Omega)$ is a consequence of
(\ref{1.4}) and the usual version of the Garsia-Rademich-Rumsey lemma
(see e.g. \cite{SV}, Theorem 2.1.3).

Consider the metric space $({\mathcal G}_p,\tilde{d}_p)$. The canonical embedding
${\mathcal H}^H \hookrightarrow {\mathcal G}_p$ is continuous. Indeed,  let $\dot{h}_i$,
$i=1,2$, belong to $L^2([0,1])$. Then for $h_i(.)=\int_0^. K(.,r) \dot{h}_i(r) dr$ and $0\leq s<t\leq 1$,
\[| (h_1)^{(1)}_{s,t} - (h_2)^{(1)}_{s,t} | \leq |t-s|^H \|\dot{h}_1-\dot{h}_2\|_2
\leq |t-s|^{\frac{1}{p}} \|h_1-h_2\|_{{\mathcal H}^H}. \]
Consequently, the preceding convergence shows that $({\mathcal G}_p, {\mathcal H}^H, P^H)$
is an abstract  Wiener space.

Let now $H\in ]\frac{1}{4}, \frac{1}{2}[$.
We follow the outline of  the proof of Lemma 3 in \cite{coutinfrizvictoir}, but
refer to the extension of the Garsia-Rademich-Rumsey lemma stated in the Lemma \ref{GRR}.

Fix $p\in ]2,4[$   such that $pH>1$. We shall prove that there exists
 $\theta >0$ such that for every $q\in [1,\infty[$,
\begin{equation}\label{eqGRR}
E\left(\left| \tilde{d}_p({\bf W,W(m)})\right|^q\right) \leq C_q 2^{-m\theta q}.
\end{equation}
Indeed, for a fixed $q\in [1,\infty[$, let $M>q$ and $N=2M$  satisfy
 $N>\frac{p}{2(Hp-1)}$. Let $\alpha,\beta>0$ defined by
$\alpha = \frac{2}{p}+\frac{1}{M}$, $\beta = \frac{1}{p}+\frac{1}{N}$.

By virtue of (\ref{1.3}) and (\ref{2.4}), we easily check that the random variables
\begin{align*}
A_1(m):&= \int_0^1 \int_0^1 ds dt  1_{\{s<t\}} \frac{| W^{(1)}_{s,t} - W(m)^{(1)}_{s,t}|^{2N}}{|t-s|^{2N\beta}},\\
A_2(m) :&= \int_0^1 \int_0^1 ds dt 1_{\{s<t\}}  \frac{| W^{(2)}_{s,t} - W(m)^{(2)}_{s,t}|^{2M}}{|t-s|^{2M\alpha}},
\end{align*}
satisfy
\begin{equation} \label{majoGRRint}
 E\big(A_1(m)\big) \leq C 2^{-m\mu 2N}\, , \; E\big(A_2(m)\big) \leq C 2^{-m\mu 2M},
\end{equation}
for some $\mu>0$.

Furthermore, the hypercontractivity property and the inequality (\ref{a16})  imply that
for $0\leq s<t\leq 1$ and $q\in [1,\infty[$,
$$\sup_m \big( \|  W^{(1)}_{s,t}\|_q + \|W(m)^{(1)}_{s,t}\|_q\big) \leq C\, |t-s|^{H}.$$
This yields
\begin{equation} \label{majoGRRintt}
 \sup_m E\big(\eta(m) \big) \leq C\, ,
\end{equation}
where
\begin{equation*}
 \eta(m):=\int_0^1 \int_0^1 ds dt 1_{\{s<t\}}  \frac{| W^{(1)}_{s,t}|^{2N} +
 | W(m)^{(1)}_{s,t}|^{2N}}{|t-s|^{2N\beta}.}
 \end{equation*}
By Lemma \ref{GRR}, we  deduce that for any  $0\leq s<t\leq 1$,
\begin{align}
 |W^{(1)}_{s,t} - W(m)^{(1)}_{s,t}|
\leq & C\, A_1(m)^{\frac{1}{2N}}\, |t-s|^{\frac{1}{p}}, \label{majoGRR1}\\
|W^{(2)}_{s,t} - W(m)^{(2)}_{s,t}|
\leq &C\,\left[ A_2(m)^{\frac{1}{2M}} +  A_1(m)^{\frac{1}{2N}}\,    \eta(m)^{\frac{1}{2N}} \right]
\, |t-s|^{\frac{2}{p}} . \label{majoGRR2}
\end{align}
Finally, Schwarz's and H\"older's inequalities together with
 (\ref{majoGRRint})-(\ref{majoGRR2}) conclude the proof of the theorem.
\hfill $\square$
\section{Appendix}

Let $W^H= (W_t^H, t\in[0,1])$ be a $d$--dimensional fractional Brownian motion with Hurst parameter
 $H\in ]0,\half[\cup]\half,1[$ and integral representation given in (\ref{a1}).

Assume $H\in]\half,1[$; by computing the integral of the right hand-side of (\ref{a3}),
 we obtain the following expression for the kernel $K^H$ defined in (\ref{a2}):
\beq
\label{a4}
K^H(t,s) = c_H\,\left(H-\half\right)s^{H-\half}\,F_2\left(\frac{t}{s}\right),
\eeq
where for $z>1$,
\beq
\label{a5}
F_2(z) = \int_0^{z-1} u^{H-\frac{3}{2}}(u+1)^{H-\frac{1}{2}}\, du.
\eeq
From (\ref{a2}), it follows that
\beq
\label{a6}
\frac{\partial K^H}{\partial t}(t,s) = c_H\,\left(H-\half\right)\,
\left(\frac{s}{t}\right)^{\half-H}(t-s)^{H-\frac{3}{2}}.
\eeq
holds  for any $H\in ]0,\half[\cup]\half,1[$ and $0<s<t<1$. Consequently,
for $H\in ]0,\half[,$
\beq
\label{a7}
\left\vert \frac{\partial K^H}{\partial t}(t,s)\right\vert \le C |t-s|^{H-\frac{3}{2}}.
\eeq
The next Lemma collects some technical estimates on the kernel $K^H(t,s)$.
\begin{lemma}
\label{la1}
Let  $0<s<t<1$.

(1) Assume $H\in ]0,\half[$. Then,
\begin{align}
&\vert K^H(t,s)\vert \le C\left(s^{H-\half}\,\1_{]0,\thalf[}(s) + (t-s)^{H-\half}\,\1_{[\thalf,t[}(s)\right),\label{a8}\\
&\left\vert \frac{\partial K^H}{\partial s}(t,s)\right\vert \le C\left(s^{H-\onehalf}\,\1_{]0,\thalf[}(s) +
(t-s)^{H-\onehalf}\,\1_{[\thalf,t[}(s)\right),\label{a9}\\
&\left\vert \frac{\partial^2 K^H}{\partial t \partial s}(t,s)\right\vert \le C (t-s)^{H-\onehalf}
\left(s^{-1}\,\1_{]0,\thalf[}(s)  + (t-s)^{-1}\,\1_{[\thalf,t[}(s)\right).\label{a10}
\end{align}

(2) For $H\in ]\half,1[$,
\begin{align}
&\vert K^H(t,s)\vert \le C \left((t- s)^{H-\half}\,\1_{]0,\thalf[}(s) +
 s^{H-\half}\,\1_{[\thalf,t[}(s)\right)\label{a11},\\
&\left\vert \frac{\partial K^H}{\partial s}(t,s)\right\vert \le C (t-s)^{2H-1}\left (s^{-(H+\frac{1}{2})}\,
 \1_{]0,\thalf[}(s) +(t-s)^{-(H+\half)}\,\1_{[\thalf,t[}(s)\right).\label{a12}
\end{align}
\end{lemma}
\begin{proof}  Assume first $H\in]0,\half[$. It is easy to check that, for any $u>0$,
\beqn
0<1-(u+1)^{H-\half}\le\left(\left(\half-H\right)u\right)\wedge 1.
\eeqn
Hence, for $0<s<t$, $0<u<\frac{t}{s}-1$,
\begin{align}
u^{H-\onehalf}\left(1-(u+1)^{H-\half}\right) &\le C u^{H-\half}\,\1_{]0,1\wedge\left(\frac{t}{s}-1\right)[}(u)\nonumber\\
&\quad+ C u^{H-\onehalf}\,\1_{]1\wedge\left(\frac{t}{s}-1\right),\frac{t}{s}-1[}(u).\label{a13}
\end{align}
Thus, from (\ref{a3}), (\ref{a13}), it follows that
\beqn
\left\vert F_1\left(\frac{t}{s}\right)\right\vert \le C\int_0^{\frac{t}{s}-1} u^{H-\frac{1}{2}} du \le C,
\eeqn
for $\frac{t}{2}\le s<t$, while for $0<s<\frac{t}{2}$,
\beqn
\left\vert F_1\left(\frac{t}{s}\right)\right\vert \le C \int_0^1u^{H-\half}\,du +C
 \int_1^\infty u^{H-\onehalf}\,du \le C.
\eeqn
Consequently
\beq
\label{a14}
\sup_{0\le s<t}\left\vert F_1\left(\frac{t}{s}\right)\right\vert \le C
\eeq
and the identity (\ref{a2}) yields (\ref{a8}).

By differentiating with respect to the variable $s$ in (\ref{a2})
and using (\ref{a14}), we obtain \beqn \left\vert \frac{\partial
K^H}{\partial s}(t,s)\right\vert \le C \left(|t-s|^{H-\onehalf}
 + s^{H-\onehalf} + s^{-1} t |t-s|^{H-\frac{3}{2}}\right),
\eeqn which yields (\ref{a9}). The inequality (\ref{a10}) follows
by differentiating with respect to the variable $s$ in (\ref{a6}).
\smallskip

Suppose now  $H\in]\half,1[$. Consider the function $F_2$ given in (\ref{a5}). Clearly, if $\frac{t}{s}-1 \le 1$,
that is, if $\frac{t}{2}\le s<t$,
\beqn
\left\vert F_2\left(\frac{t}{s}\right)\right\vert \le C.
\eeqn
Assume  $\frac{t}{s}-1 > 1$. For any $u\in]1, \frac{t}{s}-1[$, $(1+u)^{H-\half}\le C u^{H-\half}$. Consequently,
\begin{align*}
\left\vert F_2\left(\frac{t}{s}\right)\right\vert &\le C\left(\int_0^1 u^{H-\onehalf} du +
 \int_1^{\frac{t}{s}-1}u^{2H-2} du\right) 
\le C\left(\frac{t}{s}\right)^{2H-1}.
\end{align*}
The previous upper bounds, together with the representation of the kernel $K^H$ given in (\ref{a4}), imply
\begin{align*}
\vert K^H(t,s)\vert&\le C \left(s^{H-\half}\left(\frac{t}{s}\right)^{2H-1}\,\1_{]0,\frac{t}{2}[}(s)
+ s^{H-\half}\,\1_{[\frac{t}{2},t[}(s)\right)\\
&\le \left( s^{H-\half}\,\1_{]0,\frac{t}{2}[}(s)+s^{-H+\half}(t-s)^{2H-1}\,\1_{]0,\frac{t}{2}[}(s)
+ s^{H-\half}\,\1_{[\frac{t}{2},t[}(s)\right)
\end{align*}
and (\ref{a11}) follows.

Differentiating with respect to the variable $s$ in (\ref{a4}) yields
\begin{align*}
\left\vert \frac{\partial K^H}{\partial s}(t,s)\right\vert &\le C
\left(s^{H-\onehalf}F_2\left(\frac{t}{s}\right)+s^{H-\half}\frac{t}{s^2}\left(\frac{t}{s}-1\right)^{H-\onehalf}
\left(\frac{t}{s}\right)^{H-\half}\right)\\
&\le C\Big(s^{H-\onehalf}\left(\frac{t}{s}\right)^{2H-1} \1_{]0,\frac{t}{2}[}(s)
+ s^{-(H+\half)}t^{H+\half}(t-s)^{H-\onehalf}\\
&\quad +s^{H-\onehalf}\,\1_{[\frac{t}{2},t[}(s)\Big),
\end{align*}
where in the last inequality we have applied the upper bounds for
$F_2$ obtained before. Replacing in the last expression $t^{2H-1}$
by $C(s^{2H-1}+(t-s)^{2H-1})$ and $t^{H+\half}$ by
$C(s^{H+\half}+(t-s)^{H+\half})$, respectively, yields \beq
\label{a15} \left\vert \frac{\partial K^H}{\partial
s}(t,s)\right\vert \le C \left(s^{H-\onehalf}+(t-s)^{H-\onehalf}
+s^{-(H+\half)}(t-s)^{2H-1}\right). \eeq If  $0<s<\frac{t}{2}$
then, $s< t-s$ and $(t-s)^{H-\onehalf}< s^{H-\onehalf}<
s^{-(H+\half)}(t-s)^{2H-1}$, while for $\frac{t}{2}\le s<t$, the
previous inequalities are reversed accordingly. Hence (\ref{a12})
clearly follows from (\ref{a15}).
\end{proof} 
\smallskip

We introduce the notation
\begin{equation}\label{numero}
\underline{t}_m=[2^m t]2^{-m}\quad {\rm  and}\quad \overline{t}_m= \underline{t}_m + 2^{-m},
\end{equation}
for any $m\in\mathbb{N}$. Notice that,  $K_m^H$ given in (\ref{km}) satisfies $K_m^H(t,s)= 0$ if $s\ge \overline{t}_m$.
\smallskip

In the next result, we give a bound for the approximation in  quadratic mean of the kernel $K^H$ by its projection
$K_m^H$.

\begin{lemma}
\label{la2}
(1)  Let $H\in]0,\half[\cup]\half,1[$. There exists a positive constant $C$ such that for any $0<s<t\le 1$,
\begin{equation}
\sup_{m\ge 1}\int_0^1 \left(\left\vert K_m^H(t,u)-K_m^H(s,u)\right\vert^2
+\left\vert K^H(t,u)-K^H(s,u)\right\vert^2\right)\,du \le C|t-s|^{2H}.\label{a16}
\end{equation}

(2) For $H\in]0,\half[$,
\beq
\int_0^1 \left\vert K^H(t,u)-K_m^H(t,u)\right\vert^2\,du \le C\left(t\wedge2^{-m}\right)^{2H}.
\label{a17}
\eeq

(3) For $H\in]\half,1[$ and any $\lambda\in ]0, \frac{1}{2H+1}[$,
\beq
\int_0^1 \left\vert K^H(t,u)-K_m^H(t,u)\right\vert^2\,du \le C\,2^{-2m\lambda} t^{2(H-\lambda)}.
\label{a18}
\eeq
\end{lemma}
\begin{proof} The operator $\pi_m$ is a contraction on $L^2[0,1]$. Thus,
\begin{align*}
&\sup_{m\ge 1}\int_0^1 \left(\left\vert K_m^H(t,u)-K_m^H(s,u)\right\vert^2
+\left\vert K^H(t,u)-K^H(s,u)\right\vert^2\right)\,du\\
&\qquad \le 2 \int_0^1\left\vert K^H(t,u)-K^H(s,u)\right\vert^2\,du 
=2 E\left(|W_t^H-W_s^H|^2\right) =2 |t-s|^{2H},
\end{align*}
proving(\ref{a16}).

By the same argument,
\beq
\label{a19}
\int_0^1 \left\vert K^H(t,u)-K_m^H(t,u)\right\vert^2\,du\le 4\int_0^1 \vert K^H(t,u)\vert^2\,du = 4\,t^{2H}.
\eeq
Therefore (\ref{a17}) holds for $t\le C\,2^{-m}$.

Fix $t\in\Delta_I^m$ with $I>7$. We assume first $H\in]0,\half[$.  Consider the decomposition
\beq
\label{a191}
\int_0^1 \left\vert K^H(t,u)-K_m^H(t,u)\right\vert^2\,du \le C \sum_{i=1}^5 T_i(t),
\eeq
with
\begin{align*}
T_1(t)&=\int_ {0}^{t_2^m} \left\vert K^H(t,u)-K_m^H(t,u)\right\vert^2 du,\\
T_2(t)&= \int_{t_{I-3}^m}^{t_I^m}\left\vert K^H(t,u)-K_m^H(t,u)\right\vert^2 du,\\
T_3(t)&=\sum_{k=3}^{[2^{m-1}t]} \int_{\Delta_k^m} \left\vert K^H(t,u)-K_m^H(t,u)\right\vert^2 du,\\
T_4(t)&=\sum_{k=[2^{m-1}t]+2}^{I-3} \int_{\Delta_k^m} \left\vert K^H(t,u)-K_m^H(t,u)\right\vert^2 du,\\
T_5(t)&= \int_{\Delta_{[2^{m-1}t]+1}^m} \left\vert K^H(t,u)-K_m^H(t,u)\right\vert^2 du.
\end{align*}
 Schwarz's inequality and (\ref{a8}) imply
 \beqn
T_1(t)\le 4\int_{0}^{t_2^m} |K^H(t,u)|^2\, du\le C\,\int_0^{t_2^m} u^{2H-1}\, du = C\,2^{-2mH}.
\eeqn
Similarly,
 \beqn
 T_2(t)\le 4\int_ {t_{I-3}^m}^{t_I^m}\, |K^H(t,u)|^2\, du\le C\,\int_{t_{I-3}^m}^t|t-u|^{2H-1}\, du = C\,2^{-2mH}.
 \eeqn
 Let $\lambda\in]H,1[$ and $k=3,\dots,[2^{m-1}t]$, which implies $\Delta_k^m\subset]0,\thalf[$.
 By Schwarz's  inequality, the mean value theorem and (\ref{a8}), (\ref{a9}), we obtain
 \begin{align*}
& \int_{\Delta_k^m} \left\vert K^H(t,u)-K_m^H(t,u)\right\vert^2\,du \le
  2^m\int_{\Delta_k^m} du\,\int_{\Delta_k^m}\,dv
 \left\vert K^H(t,u)- K^H(t,v)\right\vert^2 \\
&\qquad \le  2^m\int_{\Delta_k^m} du\!\int_{\Delta_k^m}\!dv \left\vert K^H(t,u)- K^H(t,v)\right\vert^{2\lambda}
\left\vert |K^H(t,u)|+|K^H(t,v)|\right\vert^{2(1-\lambda)}\\
&\qquad \le C 2^{-m(2\lambda-1)}\,\int_{\Delta_k^m}du\,\int_{\Delta_k^m}\,dv\,\left((u\wedge v)^{2H-1-2\lambda}
\right).
\end{align*}
For $u,v\in \Delta_k^m$, $u\wedge v\ge u-2^{-m}$;
thus,
\[ T_3(t) 
\le C 2^{-2m\lambda }\int_{t_2^m}^{t_{[2^{m-1}t]}^m} du\, (u-2^{-m})^{2H-1-2\lambda}
\le C 2^{-2mH}. \]
Fix now $k=[2^{m-1}t]+2,\dots,I-3$, so that $\Delta_k^m\subset[\thalf, t[$. In this case
\begin{align*}
& \int_{\Delta_k^m} \left\vert K^H(t,u)-K_m^H(t,u)\right\vert^2\,du \le C 2^{-m(2\lambda-1)}\\
&\qquad\times\int_{\Delta_k^m}du\,\int_{\Delta_k^m}\,dv \left(t-(u\vee v)\right)^{2H-1-2\lambda}.
\end{align*}
Since for $u,v\in \Delta_k^m$,   $t-(u\vee v)\ge t-u-2^{-m}\ge t_{I-2}^m-u$, the previous estimate implies
\[ T_4(t) \le C 2^{-2m\lambda }\int_{t_{[2^{m-1}t]}^m}^{t_{I-3}^m}   du\,(t_{I-2}^m-u)^{2H-1-2\lambda}
\le C 2^{-2mH}. \]
We study the term $T_5(t)$ using the same method as for $T_3(t)$, $T_4(t)$, as follows:
\begin{align*}
T_5(t)&\le 2^m \int_{\Delta_{[2^{m-1}t]+1}^m}\,du \int_{\Delta_{[2^{m-1}t]+1}^m}\,dv
\left\vert K^H(t,u)-K^H(t,v)\right\vert^2\\
&\le C 2^{-m(2\lambda-1)} \int_{\Delta_{[2^{m-1}t]+1}^m}\!\!\!!du
 \int_{\Delta_{[2^{m-1}t]+1}^m}\!\!\! dv
\left((u\wedge v)^{H-\onehalf}+(t-(u\vee v)^{H-\onehalf}\right)^{2\lambda}\\
&\quad \times \left((u\wedge v)^{H-\half} + (t-(u\vee v)^{H-\half}\right)^{2(1-\lambda)}.
\end{align*}
For $u,v\in \Delta_{[2^{m-1}t]+1}^m$, $u\wedge v>\frac{t}{2}-2^{-m}$, $u\vee v<\thalf+2^{-m}$ and
$t-(u\vee v)>\frac{t}{2}-2^{-m}$. Thus, the last integral is bounded by
\beqn
 \int_{\Delta_{[2^{m-1}t]+1}^m}\,du \int_{\Delta_{[2^{m-1}t]+1}^m}\,dv \left(\frac{t}{2}-2^{-m}\right)^{2H-1-2\lambda}.
 \eeqn
Moreover, since we are assuming that $t\in\Delta_I^m$, with $I>7$, $\frac{t}{2}-2^{-m}\geq   2^{-m+1}$.
Thus, we finally obtain for $\lambda=\frac{1}{2}$,
\beqn
T_5(t)\le C 2^{-2mH}.
\eeqn
Then (\ref{a17}) follows from the upper bounds obtained so far for $T_i(t)$, $i=1,\dots,5$.
 \smallskip

 Notice that we have also proved that for $H\in ]0,\frac{1}{2}[$,
  \beq
 \label{forlemma}
 \sum_{k=3}^{I-3} 2^m\int_{\Delta_k^m} du \int_{\Delta_k^m} dv \vert K^H(t,u)-K^H(t,v)\vert^2 \le C 2^{-2mH}.
 \eeq

Assume now $H\in]\half,1[$ and fix $\lambda\in]0,\frac{1}{2H+1}[$,
so that $H-\lambda>0$. Since the inequality (\ref{a19}) holds for
any $H\in]0,\half[\cap]\half,1[$, (\ref{a18}) holds for  any $t\le
C2^{-m}$. Let now  $t\in\Delta_I^m$, with 
$I>7$. We apply a similar method as we used in the case
$H\in]0,\half[$, using the decomposition (\ref{a191}).
In fact, 
owing to (\ref{a11}),
\begin{align*}
T_1(t)&\le C \int_0^{t_2^m} (t-u)^{2H-1}\, du \le C 2^{-m} t^{2H-1},\\
T_2(t)&\le C 
\int_{t_{I-3}^m}^{t_I^m} u^{2H-1}\, du \le C 2^{-m} t^{2H-1}.
\end{align*}
Fix $k=3,\dots,[2^{m-1}t]$. Schwarz's inequality, along with the mean value theorem and (\ref{a11}), (\ref{a12}),
imply
\begin{align*}
& \int_{\Delta_k^m} \left\vert K^H(t,u)-K_m^H(t,u)\right\vert^2\,du \le  2^m\int_{\Delta_k^m} du\,\int_{\Delta_k^m}\,dv
\left\vert K^H(t,u)- K^H(t,v)\right\vert^{2\lambda}\\
&\qquad \qquad\times\left\vert |K^H(t,u)|+|K^H(t,v)|\right\vert^{2(1-\lambda)}\\
&\quad\le C 2^{-m(2\lambda-1)}\int_{\Delta_k^m} du\,\int_{\Delta_k^m}\,dv
\left((t-(u\wedge v)\right)^{(\lambda+1)(2H-1)}(u\wedge v)^{-\lambda(2H+1)}\\
&\quad\le C 2^{-2m\lambda}\, t^{(\lambda+1)(2H-1)}\int_{\Delta_k^m} du\,(u-2^{-m})^{-\lambda(2H+1)}.
\end{align*}
Since $\lambda<\frac{1}{2H+1}$, we have
\beqn T_3(t)\le C
2^{-2m\lambda}\, t^{2(H-\lambda)}.
\eeqn
Let now
$k=[2^{m-1}t]+2, \dots,I-3$. With similar arguments as before, we
deduce
\begin{align*}
& \int_{\Delta_k^m} \left\vert K^H(t,u)-K_m^H(t,u)\right\vert^2\,du \le  2^m\int_{\Delta_k^m} du\,\int_{\Delta_k^m}\,dv
\left\vert K^H(t,u)- K^H(t,v)\right\vert^{2\lambda}\\
&\qquad \qquad\times\left\vert |K^H(t,u)|+|K^H(t,v)| \right\vert^{2(1-\lambda)}\\
&\quad\le C 2^{-m(2\lambda-1)}\int_{\Delta_k^m} du\,\int_{\Delta_k^m}\,dv (t-(u\vee v))^{\lambda(2H-3)}
(u\vee v)^{(1-\lambda)(2H-1)}\\
&\quad\le C 2^{-2m\lambda}\, t^{(1-\lambda)(2H-1)}\, \int_{\Delta_k^m} du\, (t-u-2^{-m})^{\lambda(2H-3)}.
\end{align*}
For  $\lambda<\frac{1}{2H+1}$, $\lambda(2H-3)+1>0$. Hence,
\beqn
T_4(t)\le C 2^{-2m\lambda} t^{(1-\lambda)(2H-1)}
\int_{\thalf}^{t_{I-3}^m} (t-u-2^{-m})^{\lambda(2H-3)}\le C
2^{-2m\lambda}\, t^{2(H-\lambda)}. \eeqn Finally, we study the
contribution of $T_5(t)$ as follows.
\begin{align*}
T_5(t)&\le 2^{m }\int_{\Delta_{[2^{m-1}t]+1}^m}\,du
\int_{\Delta_{[2^{m-1}t]+1}^m}\,dv
\left\vert K^H(t,u)-K^H(t,v)\right\vert^2\\
&\le C2^{-m(2\lambda-1)}\int_{\Delta_{[2^{m-1}t]+1}^m}\,du \int_{\Delta_{[2^{m-1}t]+1}^m}\,dv
\Big( (t-(u\wedge v ))^{2H-1}\\
&\qquad \times\big( (u\wedge v)^{-(H+\half)}+(t-(u\vee v))^{-(H+\half)}\big)\Big)^{2\lambda}\\
&\qquad \times\Big(\big(t-(u\wedge v)\big)^{H-\half}+(u\vee v)^{H-\half} \Big)^{2(1-\lambda)}.
\end{align*}
For $u,v\in \Delta_{[2^{m-1}t]+1}^m$, $u\wedge v> C_{1}t$, $u\vee v< C_{2}t$, $t-(u\wedge v)< C_{3}t$
 and  $t-(u\vee v)> C_{4}t$. Thus,
\beqn
T_5(t)\le C\,  2^{-m(2\lambda-1)} 2^{-2m}
t^{2(H-\lambda)-1}\le C 2^{-2m\lambda}\, t^{2(H-\lambda)} \eeqn
The estimates obtained so far imply (\ref{a18}).
\end{proof} 
\smallskip

In the next Lemma we prove a simple extension of a well-known integration formula for bounded variation functions.
\begin{lemma}
\label{la6}
For any  $h\in {\IH}$, $t\ge 0$,
\beq
\label{a26}
\int_0^t h(u) h(du)= \frac{h^2(t)}{2},
\eeq
where the integral is understood in the sense of Proposition 5 in \cite{mss}.
\end{lemma}
\begin{proof} Let $n\ge 1$ and let  $h(n)$ be the function obtained by linear interpolation on the $n$-th dyadic grid
of  $h$. We have proved in \cite{mss}, Theorem 9 that
\beqn
\lim_{n\to \infty} \int_0^t h(n)(u) h(n)(du) = \int_0^t h(u) h(du),
\eeqn
for any $t\ge 0$. Since (\ref{a26}) is 
true with $h$ replaced by $h(n)$, the result follows.
\end{proof} 

The following result gives  an upper bound for the $L^2$ norm of a Skorohod integral of a Gaussian process.

\begin{lemma}
\label{la7}
Let $X_t = \int_0^1 g(t,s) dB_s$, $t\in[0,1]$, with $g$ a deterministic function belonging to $L^2([0,1]^2)$. Then, the
Skorohod integral $\int_0^1 X_s dB_s$  satisfies
\beq
\label{a27}
E\left(\int_0^1 X_s dB_s\right)^2 \le C \int_0^1 ds \int_0^1 dr  \vert g(s,r)\vert^2.
\eeq
\end{lemma}
\begin{proof} The isometry property of the Skorohod integral (\cite{nualart}, Equation (1.48)) yields
\beqn
E\left(\int_0^1 X_s dB_s\right)^2 \le C \int_0^1 E ( X_s)^2 ds + \int_0^1 ds \int_0^1 dr E(\vert D_r X_s\vert^2).
\eeqn
Since $E ( X_s)^2 = \int_0^1 \vert g(s,r)\vert^2 dr$ and the Malliavin derivative $D_r X_s$ is equal to $g(s,r)$,
(\ref{a27}) follows.
\end{proof} 

We conclude this section by proving an extension of the Garsia-Rademich-Rumsey lemma used to estimate
$\tilde{d}_p(X,Y)$ when $X$ and $Y$ are geometric rough paths with roughness $p\in[2,\infty[$ (see \cite{lyons},
Definition 3.3.3).
\begin{lemma}\label{GRR}
Let $X$ and $Y$ be geometric rough paths with the same roughness $p\in [2,+\infty[$. Set $k=[p]$.
For $i=1, \cdots, k$, let  $M_i \geq 1$, $\alpha_i = \frac{i}{p}+\frac{1}{M_i}$. Suppose that
\begin{align}
\int_0^1 \int_0^1 ds dt 1_{\{s\leq t\}}\, \frac{|X^{(i)}_{s,t}|^{2M_i}+
|Y^{(i)}_{s,t}|^{2M_i}}{|t-s|^{2M_i\alpha_i}}&
\leq A_i,  \quad 1\leq i\leq k-1, \label{X+Y}\\
\int_0^1 \int_0^1 ds dt 1_{\{s\leq t\}}\,  \frac{|X^{(i)}_{s,t} -
 Y^{(i)}_{s,t}|^{2M_i}}{|t-s|^{2M_i \alpha_i}}  &\leq B_i,\,\quad  1\leq i\leq k . \label{X-Y}
\end{align}
Then, there exists a constant $C>0$ such that for any
 $0\leq s<t\leq 1$,
\begin{align}
| X^{(i)}_{s,t}|+|Y^{(i)}_{s,t} | \leq & C\,  F_i |t-s|^{\frac{i}{p}}, \quad 1\leq i\leq k-1,
\label{GRR1}\\
\Big| X^{(i)}_{s,t}-Y^{(i)}_{s,t} \Big| \leq & C\,  G_i |t-s|^{\frac{i}{p}}, \quad 1\leq i\leq k. \label{GRR2}
\end{align}
where   $F_i$ and $G_i$ are defined recursively by
\begin{align}
F_i=&A_i^{\frac{1}{2 M_i}} + \sum_{j=1}^{i-1} F_j\, F_{i-j}, \quad 1\leq i\leq k-1, \label{F}\\
G_i=& B_i^{\frac{1}{2 M_i}} + \sum_{j=1}^{i-1} G_j\, F_{i-j}, \quad 1\leq i\leq k. \label{G}
\end{align}
\end{lemma}
\smallskip

\noindent{\bf Remark:} For rough paths $X$, $Y$ of roughness $p\in[1,\infty[$,
 $X^{(1)}_{s,t}-X^{(1)}_{s,t}= (X-Y)^{(1)}_{s,t}$. The usual version of the Garsia-Rademich-Rumsey
lemma yields the following. If \beqn \int_0^1 \int_0^1 ds dt
1_{\{s\leq t\}}\, \frac{|X^{(1)}_{s,t}-
Y^{(1)}_{s,t}|^{2M_1}}{|t-s|^{2M_1\alpha_1}} \leq B_1, \eeqn then
$| X^{(1)}_{s,t}-Y^{(1)}_{s,t} | \leq  C\,  B_1^{\frac{1}{2M_1}}
|t-s|^{\frac{1}{p}}$. Similarly, if \beqn \int_0^1 \int_0^1 ds dt
1_{\{s\leq t\}}\, \frac{|X^{(1)}_{s,t}|^{2M_1}+
|Y^{(1)}_{s,t}|^{2M_1}}{|t-s|^{2M_1\alpha_1}} \leq A_1, \eeqn then
$| X^{(1)}_{s,t}|+|Y^{(1)}_{s,t} | \leq  C\,  A_1^{\frac{1}{2M_1}}
|t-s|^{\frac{1}{p}}.$
\medskip

\noindent{\it Proof of Lemma \ref{GRR}:} Throughout the proof, the constants $F_i$, $1\leq i\leq k-1$ and
$G_i$, $1\leq i\leq k$ are defined by (\ref{F}), (\ref{G}), respectively.
We introduce the following assumption:

\noindent {\bf ($H_i$)}
\begin{align*}
&\int_0^1 \int_0^1 ds dt 1_{\{s\leq t\}}\,  \frac{|X^{(i)}_{s,t} -
 Y^{(i)}_{s,t}|^{2M_i}}{|t-s|^{2M_i \alpha_i}}  \leq B_i,\\
&| X^{(j)}_{s,t}|+|Y^{(j)}_{s,t} | \leq C\,  F_j |t-s|^{\frac{j}{p}}, \quad 1\leq j\leq i-1, \\
&| X^{(j)}_{s,t}-Y^{(j)}_{s,t}| \leq C\,  G_j |t-s|^{\frac{j}{p}}, \quad 1\leq j\leq i-1,
\end{align*}
$i\in\{2, \dots, k\}$, and we prove that {\bf ($H_i$)} implies
\beq
\label{i}
| X^{(i)}_{s,t}-Y^{(i)}_{s,t}| \leq C\,  G_i |t-s|^{\frac{i}{p}}.
\eeq
For this,
we use an argument similar to the proof of Theorem 2.1.3 in \cite{SV}.

Indeed, for every $t\in [0,1]$, set
\[ I(t)=\int_0^t  \frac{|X^{(i)}_{s,t} - Y^{(i)}_{s,t}|^{2M_i}}{|t-s|^{2M_i\alpha_i}}\, ds\, ,
\; J(t)= \int_t^1  \frac{|X^{(i)}_{t,u} - Y^{(i)}_{t,u}|^{2M_i}}{|u-t|^{2M_i\alpha_i}}\, du\, .\]
Then $\int_0^1 I(t)\, dt = \int_0^1 J(t)\, dt \leq  B_i$ and there exists $t_0>0$ such that
$I(t_0)+ J(t_0) \leq 2\, A_i$. We construct by induction a decreasing sequence $(t_n,\, n\geq 0)$ such that
$\lim_n t_n=0$ and an increasing sequence $(s_n,\, n \geq 0)$ such that $s_0=t_0$, $\lim_n s_n=1$,
 and such that there exists $C>0$ such that for every $n\geq 1$,
\begin{align}
\left| X^{(i)}_{t_n,t_0} - Y^{(i)}_{t_n,t_0}\right| \leq & C\int_0^1
\left| 8\, B_i\right|^{\frac{1}{2M_i}}\, u^{\frac{i}{p} -1}\, du + C\, \sum_{j=1}^{i-1} F_j\, G_{i-j}
\, ,\label{majot}\\
\left| X^{(i)}_{s_0,s_n} - Y^{(i)}_{s_0,s_n}\right| \leq & C\int_0^1
\left| 8\, B_i\right|^{\frac{1}{2M_i}}\, u^{\frac{i}{p} -1}\, du + C\,
\sum_{j=1}^{i-1} F_j\, G_{i-j} \, .\label{majos}
\end{align}
Then  Chen's identity  implies as  $n\rightarrow +\infty$,
\begin{align} \label{01}
|X^{(i)}_{0,1}&-Y^{(i)}_{0,1}| \leq  |X^{(i)}_{0,t_0}-Y^{(i)}_{0,t_0}| + |X^{(i)}_{t_0,1}-Y^{(i)}_{t_0,1}|
 \nonumber \\
& \quad  + \sum_{j=1}^{i-1} \left( |X^{(j)}_{0,t_0} - Y^{(j)}_{0,t_0}| |X^{(i-j)}_{t_0,1}|
+ | Y^{(j)}_{0,t_0}| |X^{(i-j)}_{t_0,1}-Y^{(i-j)}_{t_0,1}|\right) .
\end{align}
With the hypothesis {\bf ($H_i$)}, we obtain (\ref{i}) with $s=0$ and $t=1$.

To construct $(t_n)$,  we
 suppose that $t_{n-1}$ has been chosen. Let $d_{n-1}$ be defined by $d_{n-1}^{\alpha_i}=\frac{1}{2}\,
t_{n-1}^{\alpha_i}$. Then there exists $t_n\in ]0,d_{n-1}[$ such that
\[ I(t_n) \leq \frac{4\, B_i}{d_{n-1}}\quad \mbox{\rm and}\quad
 \frac{|X^{(i)}_{t_n,t_{n-1}}-Y^{(i)}_{t_n,t_{n-1}}|^{2M_i} }{|t_{n-1}-
t_n|^{2M_i\alpha_i}} \leq \frac{2 I(t_{n-1})}{d_{n-1}}\, .\]
Indeed, the sets where each one of these inequalities may fail has Lebesgue measure less that $\frac{d_{n-1}}{2}$.
Furthermore, for every $n\geq 0$, $ 2\, d_{n+1}^{\alpha_i} = t_{n+1}^{\alpha_i} \leq d_n^{\alpha_i}
=\frac{1}{2} t_n^{\alpha_i}$ and $|t_n-t_{n+1}|^{\alpha_i} \leq t_n^{\alpha_i} = 2 \, d_n^{\alpha_i}
 \leq 4\, (d_n^{\alpha_i} - d_{n+1}^{\alpha_i} )$.
Hence there exists $a\in ]0,1[$ such that $t_{n+1}\leq a\, t_n $, so that $\lim_n t_n=0$ and
more precisely,
\begin{equation}\label{geometric}
t_n \leq a^n\, t_0,
\end{equation}
while for any $n\geq 1$,
\begin{align}\label{2tn}
|X^{(i)}_{t_{n+1},t_n} - Y^{(i)}_{t_{n+1},t_n}|&\leq  |2\, I(t_n)|^{\frac{1}{2M_i}}
\, d_n^{-\frac{1}{2M_i}}\,
|t_n-t_{n+1}|^{\alpha_i} \nonumber\\
& \leq |8\, B_i|^{\frac{1}{2M_i}} |d_n\, d_{n-1}|^{-\frac{1}{2M_i}} \, 4\,
|d_n^{\alpha_i} - d_{n+1}^{\alpha_i}| \nonumber\\
&\leq  4\, \alpha_i\, \int_{d_{n+1}}^{d_n} |8\, B_i|^{\frac{1}{2M_i}} \; u^{-\frac{1}{M_i}+\alpha_i-1}\, du.
\end{align}
Let $b=a^{\frac{1}{p}} <1$;  Chen's identity,  {\bf ($H_i$)} and (\ref{2tn}) imply that for any
 $n \geq 1$,
\begin{align*}
& \left| X^{(i)}_{t_{n+1},t_0} - Y^{(i)}_{t_{n+1},t_0}\right| \leq
\left| X^{(i)}_{t_n,t_0} - Y^{(i)}_{t_n,t_0}\right|
+ \left| X^{(i)}_{t_{n+1},t_n} - Y^{(i)}_{t_{n+1},t_n}\right|\\
&\quad +   \sum_{j=1}^{i-1}
 \left( |X^{(j)}_{t_{n+1},t_n} - Y^{(j)}_{t_{n+1},t_n}|\, |X^{(i-j)}_{t_n,t_0}| +
|Y^{(j)}_{t_{n+1},t_n}|\,  |X^{(i-j)}_{t_n,t_0}- Y^{(i-j)}_{t_n,t_0}| \right) \\
& \leq \left| X^{(i)}_{t_n,t_0} - Y^{(i)}_{t_n,t_0}\right| +  C \int_{d_{n+1}}^{d_n}
 |8\, B_i|^{\frac{1}{2M_i}}\, u^{\frac{i}{p} -1}\, du \\
& \quad  + C \sum_{j=1}^{i-1} \left( G_j F_{i-j} + F_j G_{i-j}\right) |t_n-t_{n+1}|^{\frac{j}{p}}
|t_0-t_n|^{\frac{i-j}{p}} 
\end{align*}
Since $\sup_{1\leq j\leq i-1}|t_n-t_{n+1}|^{\frac{j}{p}} \leq t_n^{\frac{1}{p}}\leq  C b^n<1$,
 an easy induction on $n$  implies that for any $n\geq 1$,
\[ \left| X^{(i)}_{t_{n},t_0} - Y^{(i)}_{t_{n},t_0}\right| \leq
C \int_0^1 |8\, B_i|^{\frac{1}{2M_i}}\, u^{\frac{i}{p}-1}\, du  +
C \left( \sum_{j=1}^{i-1}  G_j F_{i-j}\right) \left( \sum_{l=0}^{n-2} b^l\right),\]
which implies (\ref{majot}). To prove (\ref{majos}), we proceed in a similar way, exchanging the endpoints
of the interval $[0,1]$.  Recall that $s_0=t_0$; suppose
 that $s_{n-1}$ has been defined and let  $\delta_{n-1}$ be such that $|1-\delta_{n-1}|^{\alpha_i}
= \frac{1}{2} |1-s_{n-1}|^{\alpha_i}$. There exists $s_n\in ]\delta_{n-1},1[$ such that
\[ J(s_n)\leq \frac{4 B_i}{1-\delta_{n-1}}\quad \mbox{\rm and}\quad
\frac{|X^{(i)}_{s_{n-1},s_n} - Y^{(i)}_{s_{n-1},s_n}|^{2M_i}}{|s_n-s_{n-1}|^{\alpha_i}} \leq
\frac{2 J(s_{n-1})}{1-\delta_{n-1}}\, .\]
Then for every $n\geq 1$, $2\, |1-\delta_{n+1}|^{\alpha_i} = |1- s_{n+1}|^{\alpha_i} \leq
 |1-\delta_n|^{\alpha_i} = \frac{1}{2} |1-t_n|^{\alpha_i}$, so that $s_n\leq \delta_n \leq s_{n+1}
 \leq \delta_{n+1}$ and for some $\bar{a} \in ]0,1[$
\begin{equation}\label{geometrics}
1- s_n \leq \bar{a}^n\, (1-t_0),
 \end{equation}
so that $\lim_n s_n=1$ and computations similar to those proving (\ref{2tn}) yield
\[
|X^{(i)}_{s_n,s_{n+1}} - Y^{(i)}_{s_n,s_{n+1}}|\leq
4\, \alpha_i\, \int_{\delta_n}^{\delta_{n+1}} |8\, B_i|^{\frac{1}{2M_i}} \; u^{-\frac{1}{M_i} + \alpha_i -1}\, du.
\]
Thus if $\bar{b}=\bar{a}^{\frac{1}{p}} <1$, Chen's identity and ($H_i$) imply
\[\left| X^{(i)}_{t_0,s_n} - Y^{(i)}_{t_0,s_n}\right| \leq
 C \int_{t_0}^{s_n} |8\, B_i|^{\frac{1}{2M_i}}\, u^{\frac{i}{p}-1}\, du  +
 C  \left( \sum_{j=1}^{i-1} F_j G_{i-j} \right) \left(\sum_{l=0}^{n-1} \bar{b}^l\right)\, ,\]
which completes the proof of (\ref{majos}) and hence that of (\ref{i}) for $s=0$, $t=1$.

To deduce (\ref{i}), for any  $s,t\in[0,1]$ with  $ s<t$, define
 $\bar{X}_u = X_{s+(t-s)u}$,
$\bar{Y}_u = Y_{s+(t-s)u}$  for $u\in [0,1]$. Then $\bar{X}$ and $\bar{Y}$ are geometric
rough paths with the same roughness $p$. Moreover,
for $0\leq u<v\leq 1$, $j=1,\cdots, k$,
 $\bar{X}^{(j)}_{u,v}=X^{(j)}_{s+(t-s)u, s+(t-s)v}$.
In fact, by a change of variables, we see that this identity is obvious for smooth rough paths
and therefore it is trivially extended to geometric rough paths.

Furthermore,
\begin{align*}
 \int_0^1 \int_0^1 du dv & 1_{\{u < v\}}\,
\frac{|\bar{X}^{(i)}_{u,v}-\bar{Y}^{(i)}_{u,v}|^{2 M_i}}{|v-u|^{2M_i\alpha_i}}\\
&\quad = (t-s)^{-2+2\alpha_i M_i}\, \int_s^t \,  \int_s^t du dv 1_{\{u < v\}}\,
 \frac{|X^{(i)}_{u,v}-Y^{(i)}_{u,v}|^{2M_i}}{|v-u|^{2M_i \alpha_i}}\\
&\quad \leq  (t-s)^{-2+2\alpha_i M_i}\, B_i = (t-s)^{2 M_i\frac{i}{p}} B_i\, .
\end{align*}

Hence, if the pair  $(X, Y)$ satisfies ($H_i$) then  $(\bar{X},\bar{Y})$ satisfies
 a similar property with constants $\bar{A}_j = (t-s)^{2M_j\frac{j}{p}}\, A_j$, $\bar{F}_j =
 |t-s|^{\frac{j}{p}} F_j$,$1\le j\le i-1$,  $\bar{B}_j= (t-s)^{2 M_j\frac{j}{p}} B_j$,
$\bar{G}_j=|t-s|^{\frac{j}{p}}$, $1\le j\le i$. This finishes the proof of (\ref{i}).

Taking in the preceding arguments first $X\equiv 0$ and then $Y\equiv 0$, we see  recursively
 that (\ref{X+Y}) implies  {\bf($H_i$)}  for any $i=1,\dots, k-1$, with $B_i=A_i$.
Hence we obtain (\ref{GRR1}). Moreover, we also see  that
{\bf($H_i$)} holds true for any $i=1,\dots, k$, whenever (\ref{X+Y}), (\ref{X-Y}) are
satisfied.
This concludes the proof.
\end{proof}
\medskip

\noindent {\bf Acknowledgments}
The first named author whishes to thank the Centre de Recerca Matem\`atica in Bellaterra
and the Universitat de Barcelona for their support and hospitality in the fall of 2004,
 when discussions on the content of this paper started.


\begin{thebibliography}{99}
\bibitem{amn} E. Al\`os, O. Mazet, D. Nualart: Stochastic calculus with respect to Gaussian processes.
{\it The Annals of Probab.} Vol. 29, No. 2, 766-801 (2001).
\bibitem{coutinquian} L. Coutin, Z. Qian: Stochastic analysis, rough path analysis and fractional
Brownian motions. {\it Probab. Theory Relat. Fields} 122, 108-140 (2002).
\bibitem{coutinfrizvictoir} L. Coutin, P. Friz, N. Victoir: Good rough path sequences and applications
to anticipating and fractional stochastic calculus.
arXiv:math.PR/0501197, January 2005.
\bibitem{du}  L. Decreusefond, S. \"Ust\"unel: Stochastic analysis of the fractional Brownian
motion. {\it Potential Analysis} 10, 177-214 (1999).
\bibitem{frizvictoir1} P. Friz, N. Victoir: A note on the notion of geometric rough paths. arXiv:math.FA/0403115, March 2004.
\bibitem{lyons} T. Lyons, Z. Qian: System Control and Rough Paths. Oxford Mathematical Monographs.
Oxford Science Publiccations. Clarendon Press, Oxford 2002.
\bibitem{mss} A. Millet, M. Sanz-Sol\'e: Large deviations for rough paths of the fractional Brownian motion.
Annales de l'Institut Poincar\'e (to appear), arXiv:math.PR/04122000, December 2004.
\bibitem{nualart} D. Nualart: The Malliavin Calculus and Related Topics. Probability and its Applications.
 Springer Verlag, 1995.
\bibitem{SV} D.W. Stroock, S.R.S. Varadhan: Multidimensional Diffusion Processes. Grundlehren des
mathematischen Wissenschaften 233. Springer Verlag, 1979.
\end{thebibliography}
\end{document}